\documentclass{amsart}
\usepackage{stackengine}
\usepackage{amssymb}
\usepackage{amsfonts}
\usepackage{amsmath}
\usepackage{amstext}
\usepackage{amsthm}
\usepackage{latexsym}
\usepackage{bbm}
\usepackage{mathrsfs}  
\usepackage[mathscr]{eucal}

\usepackage{amscd}
\usepackage{eucal}
\usepackage{amsxtra}
\usepackage{amsbsy}
\usepackage{graphicx}
\usepackage{latexsym}
\usepackage{pb-diagram}
\usepackage{amsopn}
\usepackage{delarray}
\usepackage{psfrag}
\usepackage{verbatim}
\usepackage{bussproofs}
\usepackage[T1]{fontenc}
\usepackage{color,cancel}

\usepackage{appendix}
\usepackage{color}
\usepackage[dvipsnames]{xcolor}
\usepackage{epigraph}
\usepackage{lipsum}
\usepackage{hyperref} 
\hypersetup{pdfpagemode={UseOutlines},
bookmarksopen=true,bookmarksopenlevel=0,hypertexnames=false,colorlinks=true,citecolor=black,linkcolor=black,urlcolor=black,pdfstartview={FitV},unicode,breaklinks=true}
\usepackage{graphicx}
\usepackage{enumitem}

\usepackage{mathrsfs}

\makeatletter
\@namedef{subjclassname@2020}{%
  \textup{2020} Mathematics Subject Classification}
\makeatother

\def\cD{{\mathcal{D}}}

\def\cI{{\mathcal{I}}}
\def\cJ{{\mathcal{J}}}

\def\cP{{\mathcal{P}}}

\def\cA{{\mathcal{A}}}
\def\cB{{\mathcal{B}}}

\newcommand{\conc}{^\smallfrown}

\newcommand{\w}{\omega}

\newcommand{\finishclaim}{\hfill\ensuremath{_{Claim}\square}}

\def\cA{{\mathcal{A}}} \def\cB{{\mathcal{B}}}  \def\cD{{\mathcal{D}}}     \def\cI{{\mathcal{I}}} \def\cJ{{\mathcal{J}}}      \def\cP{{\mathcal{P}}} \def\cQ{{\mathcal{Q}}}         

                \def\bQ{{\mathbb{Q}}}

\def\P{{\mathbb{P}}}
\def\Q{{\mathbb{Q}}}

\def\dom{\mathrm{dom}}
\def\LIM{\mathrm{LIM}}

\def\DG{\mathscr{DG}}

\newcommand{\diam}{\diamondsuit}
\newcommand{\diamT}{\diamondsuit_T}
\newcommand{\up}{\mathord{\uparrow}}

\newcommand{\II}{\mathrm{II}}
\newcommand{\I}{\mathrm{I}}

\newcommand{\down}{\mathord{\downarrow}}
\newcommand{\forces}{\Vdash}

\newcommand{\wsic}{wscm}
\newcommand{\sic}{scm}
\newcommand{\wsicc}{wscm }
\newcommand{\sicc}{scm }

\newcommand{\vecD}{\overrightarrow{\mathscr{D}}}
\newcommand{\htt}{\operatorname{ht}}
\definecolor{green(html/cssgreen)}{rgb}{0.0, 0.5, 0.0}

\definecolor{britishracinggreen}{rgb}{0.0, 0.26, 0.15}

\newcommand{\vecC}{\overrightarrow{C}}

\DeclareMathOperator{\htop}{ht}

\newcommand{\Katetov}{Kat\v{e}tov }

\newtheorem{theorem}{Theorem}[section]
\newtheorem{proposition}[theorem]{Proposition}
\newtheorem{corollary}[theorem]{Corollary}
\newtheorem{lemma}[theorem]{Lemma}
\newtheorem{question}[theorem]{Question}
\newtheorem*{questionn}{Question}

\theoremstyle{definition}
\newtheorem{mydef}[theorem]{Definition}

\theoremstyle{remark}
\newtheorem{remark}[theorem]{Remark}
\newtheorem{fact}[theorem]{Fact}

\newtheorem{claim}[theorem]{Claim}
\newtheorem{subclaim}[theorem]{Subclaim}

\newtheorem{notation}[theorem]{Notation}

\theoremstyle{remark}

\begin{document}


\baselineskip=17pt


\title[Diamonds on trees]{Diamonds on trees}

\author[O. Guzm\'an]{Osvaldo Guzm\'an}
\address{Centro de Ciencias  Matem\'aticas, Universidad Nacional Aut\'onoma de M\'exico, Campus Morelia, 58089, Morelia, Michoac\'an, M\'exico.}
\email{oguzman@matmor.unam.mx}

\author[C. L\'opez-Callejas]{Carlos L\'opez-Callejas}
\address{Centro de Ciencias  Matem\'aticas, Universidad Nacional Aut\'onoma de M\'exico, Campus Morelia, 58089, Morelia, Michoac\'an, M\'exico.}
\email{carloscallejas@matmor.unam.mx}



\date{}


\begin{abstract}
We generalize the diamond principle and its variants using the notion of stationarity in trees introduced by Brodsky in \cite{brodskyarticle} and \cite{brodskyPhDthesis}. In particular, we show that if $T$ is a nonspecial $\omega_1$-tree, then $\diamT \implies \diam$, and if $T$ is a Suslin tree, then $\diamT \iff \diam$. We also prove that $\diam^*$ implies $\diamT$ (yielding the consistency of $\diamT$) and establish the consistency of $\neg\diam^* + (\forall T\text{ nonspecial }\omega_1\text{-tree }(\diamT))$. Finally, we demonstrate that it is consistent with $\diam$ that there exists a nonspecial $\omega_1$-tree with \(\neg\diamT\), introducing two forcing properties—$\sigma(S)$-closed and strategically closed in models—which are preserved under countable support iterations.
\end{abstract}

\maketitle


\section{Introduction}

In \cite{stevostationary} and \cite{stevopartitionrelations}, Stevo Todorcevic introduced a new notion of stationarity for subsets of \(\omega_1\) parameterized by an arbitrary nonspecial tree of height \(\omega_1\). He proved that such a tree \(T\) is special if and only if there exists a regressive map \(f\colon T\to T\) such that \(f^{-1}(\{t\})\) is special for every \(t\in T\). This notion was later used by Ari Brodsky in \cite{brodskyPhDthesis} and \cite{brodskyarticle} to define a notion of nonstationarity for subsets of \(T\) itself\footnote{Brodsky's motivation was to extend the well-known Balanced Baumgartner--Hajnal--Todorcevic Theorem to the setting of trees.}, rather than for subsets of its height \(\omega_1\). Thus, a subset \(X\subseteq T\) is nonstationary if and only if there is a regressive map \(f\colon X\to T\) such that \(f^{-1}(\{t\})\) is special for all \(t\in T\).

Once the notion of stationarity is in place, it becomes natural to define diamond sequences on trees. More precisely, one considers sequences of the form $(D_t \mid t \in T)$, where for each node $t \in T$, the set $D_t$ is a subset of $\{ s \in T \mid s < t \}$. These sequences are required to \textit{guess} every subset $X \subseteq T$ in a stationary way. We denote such sequences by $\diamT$. In this work, we study $\diamT$-sequences in depth and explore the relationships between these and the usual diamond principles defined on $\omega_1$. 


Our main results include the following. We show that if $T$ is an $\omega_1$-tree that either has a cofinal branch or is almost-Suslin (in particular, if $T$ is Suslin), then $\diamT$ is equivalent to the classical diamond principle $\diam$. This gives rise to the informal intuition that for Aronszajn trees, the \textit{thinner} the tree, the more $\diamT$ resembles $\diam$. On the other hand, we prove that $\diam_{2^{<\w_1}}$ holds in $\mathsf{ZFC}$. 
This gives rise to the informal intuition that the more cofinal branches a tree $T$ has, the more likely it is that $\diamT$ holds\footnote{Curiously, the complete opposite is true for another diamond principle for trees, $\diamondsuit(T)$, studied in \cite{paper66}.}.

We also prove that $\diam^*$ implies $\diamT$ for every nonspecial and well-pruned $\omega_1$-tree $T$; in particular, we prove an analogue of Kunen's theorem stating that the principles $\diam$ and $\diam^-$ are equivalent; then we show that $\diam^*$ turns out to be strictly stronger than the principle asserting that $\diamT$ holds for every nonspecial $\omega_1$-tree $T$. 

We also establish the consistency of $\diam$ together with the failure of $\diam_T$ for some nonspecial Aronszajn tree $T$. This is achieved by starting from a specific nonspecial Aronszajn tree in a model of $V=L$, and by showing that in a suitable forcing extension $\diam$ is preserved while $\diam_T$ is destroyed (without specializing $T$). Hence, for some Aronszajn trees $\diam_T$ lies strictly between $\diam$ and $\diam^*$, demonstrating that $\diam_T$ is a genuinely new principle.

We view $\diamT$ as an \emph{invariant} for studying nonspecial $\omega_1$-trees, in the sense that it behaves monotonically  between trees, that is, if there is a strictly increasing map $f:S\to T$, then $\diam_S$ implies $\diamT$ (see part (ii) of Theorem \ref{Lipschitz map and katetov}). In other words, our study of $\diamT$-sequences is not primarily aimed at finding applications, but rather at using them as invariants to better understand the trees themselves. Something that supports this paradigm for studying trees is that through $\diamT$ we can consider three types of nonspecial trees of height $\omega_1$: those for which $\diamT$ holds in $\mathsf{ZFC}$ (such as $2^{<\omega_1}$ (see Theorem \ref{diam for 2^<w1 is true in ZFC})), those for which $\diamT$ is equivalent to the classical $\diam$ (such as Suslin trees; see Corollary~\ref{diam iff diamT if T is Suslin}), and those for which $\diamT$ is strictly stronger than the classical $\diam$ (such as the one exhibited in Theorem~\ref{conjectura diamante mas debil que diamante T}).

The article is organized as follows. In Section 2, we compile all the necessary definitions and results from Brodsky’s work on stationarity in trees (without proofs), and we include a few additional results that will be needed later and that are not explicitly found in his papers. In Section 3, we introduce the concept of $\diamT$ and investigate some of its basic properties; we start by studying how $\diam_S$ and $\diamT$ compare when $S$ is a subtree of $T$, then we prove that $\diamT$ imposes certain cardinal restrictions on the size of $T$ in the same way that $\diam$ implies $\mathsf{CH}$, that $\diam_{2^{<\omega_1}}$ holds in $\mathsf{ZFC}$, and that $\diamT$ implies $\diam$ for every nonspecial $\omega_1$-tree. We conclude the section by showing that $\diamT$ can be expressed within a more familiar framework, commonly used to define the diamond principle on $P_\kappa(\lambda)$. In Section 4, we prove the consistency of the statement $\forall T\, \Bigl(\text{$T$ is a nonspecial $\omega_1$-tree} \implies \diamT\Bigr)$ and we show that $\diamT\iff\diam$ when $T$ is an almost-Suslin tree. In Section 5 we construct a model of $\neg\diam^*$ where $\diamT$ holds for every $T$  nonspecial $\omega_1$-tree. In Section 6 we separate the notion of $\diamT$ from the classical diamond principle $\diam$ for some Aronszajn tree $T$. In this section we introduce the two forcing properties —namely, the $\sigma(S)$-closed and the strategically closed in models forcings—and prove iteration theorems for both. Finally we conclude with a section of open questions.

\section{Preliminaries and Brodsky's notion of stationary for trees}

We begin by recalling some classical guessing principles that will appear throughout this article, starting with Jensen’s notion of diamond:

\begin{mydef}[\cite{jensen}]
Let $S\subseteq\omega_1$ be stationary. A sequence $(D_\alpha)_{\alpha\in S}$ is a \emph{$\diam(S)$-sequence} if
\begin{enumerate}
  \item $D_\alpha\subseteq\alpha$ for every $\alpha\in S$, and
  \item for every $X\subseteq\omega_1$, the set $\{\alpha\in S \mid X\cap\alpha=D_\alpha\}$ is stationary.
\end{enumerate}
\end{mydef}

Two related principles that strengthen $\diam(S)$ are:

\begin{mydef}[\cite{jensennotes}]
Let $S\subseteq\omega_1$ be stationary. A sequence $(\mathcal D_\alpha)_{\alpha\in S}$ is a \emph{$\diam^*(S)$-sequence} if
\begin{enumerate}
  \item $\mathcal D_\alpha\subseteq\mathcal P(\alpha)$ and $|\mathcal D_\alpha|\le\omega$ for every $\alpha\in S$, and
  \item for every $X\subseteq\omega_1$, the set $\{\alpha\in S \mid X\cap\alpha\in\mathcal D_\alpha\}$ contains a club.
\end{enumerate}
\end{mydef}

\begin{mydef}[\cite{jensennotes}]
Let $S\subseteq\omega_1$ be stationary. A sequence $(\mathcal D_\alpha)_{\alpha\in S}$ is a \emph{$\diam^+(S)$-sequence} if
\begin{enumerate}
  \item $\mathcal D_\alpha\subseteq\mathcal P(\alpha)$ and $|\mathcal D_\alpha|\le\omega$ for every $\alpha\in S$, and
  \item for every $X\subseteq\omega_1$, there is a club $C\subseteq\omega_1$ such that for all $\alpha\in C\cap S$ we have $X\cap\alpha\in\mathcal D_\alpha$ and $C\cap\alpha\in\mathcal D_\alpha$.
\end{enumerate}
\end{mydef}

\noindent It is clear that $\diam^+(S)\implies\diam^*(S)$; a natural weakening of $\diam^*(S)$ is the following:

\begin{mydef}
Let $S\subseteq\omega_1$ be stationary. A sequence $(\mathcal D_\alpha)_{\alpha\in S}$ is a \emph{$\diam^-(S)$-sequence} if
\begin{enumerate}
  \item $\mathcal D_\alpha\subseteq\mathcal P(\alpha)$ and $|\mathcal D_\alpha|\le\omega$ for every $\alpha\in S$, and
  \item for every $X\subseteq\omega_1$, the set $\{\alpha\in S \mid X\cap\alpha\in\mathcal D_\alpha\}$ is stationary.
\end{enumerate}
\end{mydef}

A remarkable result of Kunen is that $\diam(S)$ and $\diam^-(S)$ are equivalent \cite[Theorem 7.14]{kunenbook}.  Consequently, $\diam^+(S)\implies\diam^*(S)\implies\diam(S)$, and none of these implications can be reversed (see \cite{devlintrees} and \cite{devlinvariations}).

Now, we present the definitions and results from Brodsky’s work that are needed for the development of this article,  for this we will fix some notation for trees and ideals. 

Recall that a tree is a partial order $(T,<_T)$ such that for every $t\in T$, the set of predecessors of $t$ is well-ordered by $<_T$. Throughout, we assume every tree has a root, which we identify—without loss of generality—with the empty set $\emptyset$.

\begin{notation}
Let $T$ be a tree.
\begin{enumerate}
  \item For each $t\in T$,
  \[
    t\down = \{\,s\in T \mid s<_T t\},
    \quad
    t\up =
    \begin{cases}
      \{\,s\in T \mid t<_T s\}, & t\neq\emptyset,\\
      T,                         & t=\emptyset.
    \end{cases}
  \]
  \item For $t\in T$, $\htop(t)$ is the \emph{height} of $t$, i.e.,\ the order type of $t\down$.
  \item For $S\subseteq T$, $\widehat S = \{\,\htop(s)\mid s\in S\}$.
  \item For a set of ordinals $C$, $T\restriction C = \{\,t\in T\mid \htop(t)\in C\}$.
  \item For an ordinal $\alpha$, 
    $T_{\alpha} = T\restriction\{\alpha\}$,
    $T_{<\alpha} = T\restriction\alpha$ and $T_{\geq\alpha}=T\setminus T_{<\alpha}$.
\end{enumerate}
\end{notation}

The \emph{height} of $T$ is the least ordinal $\alpha$ with $T_{\alpha}=\emptyset$, denoted by $\htop(T)$.  A \emph{$\kappa$-tree} is a tree of height $\kappa$ in which each level $T_{\alpha}$ has size $<\kappa$.  Thus an $\omega_{1}$-tree has height $\omega_{1}$ and each level is at most countable.

A \emph{cofinal branch} in $T$ is a chain $B\subseteq T$ meeting every level: $B\cap T_{\alpha}\neq\emptyset$ for all $\alpha<\htop(T)$.  An \emph{antichain} is a subset $A\subseteq T$ in which no two distinct elements are comparable: if $s\neq t\in A$ then neither $s<_T t$ nor $t<_T s$.  An $\omega_{1}$-tree is \emph{Aronszajn} if it has no cofinal branches, and \emph{Souslin} if it has neither cofinal branches nor uncountable antichains. For a tree $T$ of height $\kappa$, we say it is \emph{well-pruned} if for every $\alpha<\beta<\kappa$ and every $s\in T_\alpha$ there exists $t\in T_\beta$ with $s<_T t$. We say $T$ is \emph{Hausdorff} if for all distinct $s,t\in T$, if $\htop(s)$ and $\htop(t)$ are limit ordinals, then $s\down\neq t\down$. Finally, if $T$ is a tree, $X\subseteq T$ and $f:X\to T$, then $f$ is \emph{regressive} if $f(t)<_T t$ for every $t\in X\setminus\{\emptyset\}$ (see \cite[Section 1]{stevostationary}).

\begin{mydef}
Let $X$ be a set. An \emph{ideal} on $X$ is a nonempty family $\mathcal I\subseteq\mathcal P(X)$\footnote{\(\cP(X)\) denotes the power set of \(X\), that is, the family of all subsets of $X$.} such that:
\begin{enumerate}
  \item If $A,B\in\mathcal I$, then $A\cup B\in\mathcal I$.
  \item If $A\subseteq B$ and $B\in\mathcal I$, then $A\in\mathcal I$.
\end{enumerate}
\end{mydef}

 For an ideal $\cI$ on $X$ its \emph{positive} part is $
\mathcal I^+=\{A\subseteq X\mid A\notin\mathcal I\}$ while
and its \emph{dual} is \mbox{$\mathcal I^*=\{A\subseteq X\mid X\setminus A\in\mathcal I\}$}. We call an ideal $\mathcal I$ \emph{proper} if $X\notin\mathcal I$; note that $\mathcal I$ is proper exactly when $\mathcal I^*\cap\mathcal I=\emptyset$, equivalently, when $\mathcal I^*\subseteq\mathcal I^+$. If $Y\subseteq X$, then the restriction $\cI\restriction Y$ is defined by
$\cI\restriction Y = \{A\cap Y \mid A\in\cI\}$ (equivalently,
$\cI\restriction Y = \{A\in\cI \mid A\subseteq Y\}$); it is easy to see that $\cI\restriction Y$ is an ideal on $Y$ and it is proper if and only if $Y\in\cI^+$. We say that $\cI$ is \emph{$\kappa$-complete}, where $\kappa$ is a cardinal, if $\bigcup\mathcal F\in\mathcal I$ for every $\mathcal F\subseteq\mathcal I$ with $|\mathcal F|<\kappa$.

The rest of the terminology we follow is standard and largely agrees with Brodsky’s. The only difference is purely linguistic: in Brodsky’s approach, every subset of a tree is called a “subtree”, whereas we reserve that term for subsets that are closed under initial segments. Thus, if $T$ is a tree and $S \subseteq T$, we say that $S$ is a \emph{subtree} of $T$, if whenever $s, t \in T$ satisfy $s \le_T t$ and $t \in S$, then $s \in S$.

Before defining nonstationary subsets of a tree, we first recall the notion of special subsets:

\begin{mydef}
    Let $T$ be a tree of height $\omega_1$ and let $U\subseteq T$. We say that $U$ is \emph{special} if $U$ can be written as the union of at most $\omega$ antichains; equivalently, there exists a function
    \[
    f:U\longrightarrow\omega
    \]
    such that for all $t,u\in U$, $t<_T u$ implies $f(t)\neq f(u)$.
\end{mydef}

Note that the collection of special subsets of $T$ forms an ideal on $T$ that is, in fact, $\omega_1$-complete. Also, it is important to observe that a special subset $U\subseteq T$ need \emph{not} be a subtree; that is, $U$ is \emph{not} required to be closed under initial segments.  

Special sets are natural analogues of the bounded subsets of $\omega_1$ in the sense that, when $\omega_1$ is viewed as a tree, the ideal of special subtrees coincides with the ideal of bounded (i.e., at most countable) subsets.

Now, we can explicitly formulate Todorcevic's notion of nonstationarity parametrized by a tree $T$ of height $\omega_1$:

\begin{mydef}\label{def:stevo-NS}\cite{stevostationary,stevopartitionrelations}
    Let $T$ be a tree of height $\omega_1$. The ideal $NS_T$ on $\omega_1$ is defined as follows: for $X\subseteq\omega_1$ we have $X\in NS_T$ if and only if there exists a regressive function
    $f\colon T\restriction X\to T$ such that $f^{-1}(\{t\})$ is special for every $t\in T$.
\end{mydef}

And its natural generalization, due to Brodsky, which motivated our investigation of guessing principles on trees:\footnote{Actually, Todorcevic treats trees of any regular uncountable height, while Brodsky does so for any successor cardinal. Here we restrict to the case of $\omega_1$, as it is the context of this article.}

\begin{mydef}\label{def:brodsky-NS}\cite{brodskyPhDthesis,brodskyarticle}\label{definiton of non-stationary}
    Let $T$ be a tree of height $\omega_1$. The ideal $NS^T$ on $T$ is defined as follows: for $B\subseteq T$ we have $B\in NS^T$ if and only if there exists a regressive function $f\colon B\to T$ such that $f^{-1}(\{t\})$ is special for every $t\in T$. In case $B\in NS^T$, we say that $B$ is \emph{nonstationary in $T$}.
\end{mydef}

Clearly every special subset of $T$ is nonstationary in $T$. Additionally, Definition \ref{def:brodsky-NS} generalizes Definition \ref{def:stevo-NS} in the sense that, for a tree $T$ of height $\omega_1$ and $X\subseteq\omega_1$, we have
\[
X\in NS_T \iff T\restriction X\in NS^T.
\]
Note also that when $T=\omega_1$, we have $NS^T=NS_T$ and, by Fodor's Lemma \cite{fodor}, both coincide with the classical nonstationary ideal on $\omega_1$.

For many arguments that follow, it is convenient to have an equivalent description of $NS^T$ in terms of \emph{diagonal unions}\footnote{Actually, Brodsky's original formulation of the ideal $NS^T$ is in terms of diagonal unions.}. 

\begin{mydef}
    Let $T$ be a tree of height $\w_1$. For any collection of subsets of $T$ indexed by the nodes of $T$, i.e., $(A_t)_{t\in T}$, its diagonal union is defined as: 
    \[
    \nabla_{t\in T} A_t=\bigcup_{t\in T}(A_t\cap t\up).
    \]
\end{mydef}

\begin{remark}\label{caracterization of diagonal union}\cite[Lemma 37]{brodskyPhDthesis}
Let $T$ be a tree of height $\w_1$ and $(A_t)_{t\in T}$ a sequence of subsets of $T$. Then 
    $\nabla_{t\in T}A_t=\{w\in T\mid w\in A_\emptyset\cup \bigcup_{t<_T w}A_t\}$. 
\end{remark}

The following theorem collects several key facts about the ideal $NS^T$ that will be used throughout the paper.

\begin{theorem}\label{teo que junta cosas sobre NS^T}
    Let $T$ be a tree of height $\w_1$. Then:
    \begin{enumerate}[label=\textup{(\Roman*)}]
        \item For every $B\subseteq T$, $B\in NS^T$ if and only if $B$ can be written as a diagonal union of special subsets of $T$. \cite[Lemmas 42 and 45]{brodskyPhDthesis} (cf.\ \cite[p.~9]{structural})
        \item\label{ideal of non-stationary is complete} $NS^T$ is $\w_1$-complete. \cite[Lemma 46]{brodskyPhDthesis}
        \item $NS^T$ is closed under diagonal unions; namely, if $\{A_t\mid t\in T\}\subseteq NS^T$, then $\nabla_{t\in T} A_t\in NS^T$. \cite[Theorem 47]{brodskyPhDthesis}
        \item If $X\in (NS^T)^+$ and $f:X\to T$ is regressive, then there exists $Y\in (NS^T)^+\cap\cP(X)$ such that $f\restriction Y$ is constant. \cite[Theorem 47 and Corollary 44]{brodskyPhDthesis}
        \item\label{pressing down for trees} 
       $T$ is nonspecial if and only if $NS^T$ is a proper ideal if and only if $NS_T$ is a proper ideal. \cite[Theorem 49]{brodskyPhDthesis}\cite{stevostationary,stevopartitionrelations}
    \end{enumerate}
\end{theorem}

Additionally, the following theorem summarizes some of the most important results relating the ideal $NS^T$ to the usual ideal of nonstationary sets in $\omega_1$. 

\begin{theorem}\label{restrictionclubnotstationary}
    \cite[Lemma 51]{brodskyPhDthesis}\cite{stevostationary,stevopartitionrelations} Let $T$ be a tree of height $\w_1$ and $X,C\subseteq \w_1$. Then:
    \begin{enumerate}[label=\textup{(\roman*)}]
        \item If $|X|\leq\w$ the $T\restriction X\in NS^T$.
        \item\label{2 de restrictionclubnotstationary} If $X$ is a nonstationary subset of $\w_1$, then $X\in NS_T$ and therefore $T\restriction X\in NS^T$.
        \item In particular, the set of successor nodes of $T$ is a nonstationary subset of $T$.
        \item\label{4 de restrictionclubnotstationary} If $C$ is a club subset of $\w_1$, then $T\restriction C\in (NS^T)^*$ and $C\in (NS_T)^*$.
        \item\label{5 de restrictionclubnotstationary} If $T$ is a nonspecial tree and $C$ is a club subset of $\w_1$, then $T\restriction C\not\in NS^T$.
    \end{enumerate}
\end{theorem}

From now on, the results in this section are no longer (at least not in an obvious manner) found in Brodsky's work. Therefore, even though these results are straightforward, we have decided to include their proofs.

In the remainder of this section, for a tree $T$ and $S\subseteq T$, whenever $t\in S$ we denote by $(t\up)_S$ the set $(t\up)\cap S$, and analogously with $t\down_S$. Also, in the following proofs we use the fact that if $X\subseteq S$, then $X$ is special in $S$ if and only if it is special in $T$, since $S$ inherits its order from $T$.

\begin{proposition}\label{NS^S contenid en NS^T restricted to S}
    Let $T$ be a tree of height $\w_1$ and $S\subseteq T$, then $NS^S \;\subseteq\; (NS^T)\restriction S$.\footnote{Here we view $S$ with the order inherited from $T$; hence, $S$ is also a tree (although it is not a subtree of $T$) and it makes sense to consider $NS^S$ (making the convention $NS^S=\cP(S)$ in case $\htop(S)\leq\w$).}
\end{proposition}

\begin{proof}
    Let $(X_s)_{s\in S}$ be such that each $X_s$ is a special subset of $S$. Now for every $t\in T$ define $Y_t=X_t$ in the case that $t\in S$ and $Y_t=\emptyset$ otherwise; note that $Y_t\subseteq S$ for all $t\in T$. 
    \begin{claim}
        $(\nabla_{t\in T}Y_t)\cap S=\nabla_{s\in S}X_s$.
    \end{claim}
    \textit{Proof of the Claim.} 
    $(\subseteq)$ If $u\in(\nabla_{t\in T}Y_t)\cap S$ there there is $t\in T$ such that $u\in Y_t\cap (t\up)_T$, then $Y_t\neq\emptyset$ and consequently $t\in S$, then $u\in X_t\cap (t\up)_T$, but as $u\in S$ then $u\in (t\up)_S$, then $u\in X_t\cap (t\up)_S\subseteq \nabla_{s\in S}X_s$.

    $(\supseteq)$ If $u\in\nabla_{s\in S}X_s$, then there is $s\in S$ such that $u\in X_s\cap (s\up)_S$, then $u\in X_s\cap (s\up)_T$ as $(s\up)_S\subseteq (s\up)_T$, then $u\in\nabla_{t\in T}Y_t$.
    \hfill\ensuremath{_{Claim}\square}
\end{proof}

Recall that $S \subseteq T$ is a subtree of $T$ if it is closed under initial segments.

\begin{lemma}\label{non statinary for subtrees is just restriction}
    Let $T$ be a tree of height $\w_1$. If $S$ is a subtree of $T$, then $NS^S=NS^T\restriction S$.
\end{lemma}
\begin{proof}
    By Proposition \ref{NS^S contenid en NS^T restricted to S}, it is enough to prove that $NS^T\restriction S\subseteq NS^S$. Let $(X_t)_{t\in T}$ be a sequence of special subsets. For each $s\in S$, let $Y_s=X_s\cap S$. Clearly, each $Y_s$ is special, as it is a subset of a special set. 

    \begin{claim}
        $\nabla_{s\in S}Y_s=(\nabla_{t\in T}X_t)\cap S$.
    \end{claim}
    \textit{Proof of the Claim.} 
    ($\subseteq$) Let $u\in Y_s\cap (s\up)_S$ for some $s\in S$, then $u\in X_s\cap (s\up)_T$ and thus $u\in\nabla_{t\in T}X_t$.

    ($\supseteq$) Let $u\in(\nabla_{t\in T}X_t)\cap S$. Thus, by Remark \ref{caracterization of diagonal union} we have $u\in X_\emptyset$ or $\exists t\in T(t<u\wedge u\in X_t)$. If $u\in X_\emptyset$, then $u\in Y_\emptyset=X_\emptyset\cap S$. If $\exists t\in T(t<u\wedge u\in X_t)$, then $t\in S$ as $S$ is a subtree of $T$, then $u\in X_t\cap (t\up)_S=Y_t\cap (t\up)_S\subseteq\nabla_{s\in S}Y_s$.
    \hfill\ensuremath{_{Claim}\square}
\end{proof}

Note that part \ref{pressing down for trees} of Theorem \ref{teo que junta cosas sobre NS^T} together with Lemma~\ref{non statinary for subtrees is just restriction} implies the following.

\begin{corollary}\label{for subtrees special and nonstationary is the same}
    If $T$ is a tree of height $\w_1$ and $S$ is a subtree of $T$, then $S$ is special if and only if $S\in NS^T$.
\end{corollary}

\section{Diamonds on trees and its basic properties}

In light of Definition \ref{definiton of non-stationary} and the fact that it includes the case $T = \omega_1$, the following guessing principles are natural generalizations of the usual ones and represent the central objects of study in this article:

\begin{mydef}
    If $T$ is a nonspecial tree of height $\w_1$. A sequence $(D_t\mid t\in T)$  is called a $\diamondsuit_T$-sequence if:
    \begin{enumerate}[label=(\arabic*)]
        \item $D_t\subseteq t\down$ for every $t\in T$.
        \item\label{condiciondeadivinanza1} For every $X\subseteq T$ we have $\{t\in T\mid X\cap t\down=D_t\}\not\in NS^T$. 
    \end{enumerate}
\end{mydef}

We refer to condition \ref{condiciondeadivinanza1} as \textit{every set $X$ is guessed in a stationary set.} 

\begin{mydef}
    If $T$ is a nonspecial tree of height $\w_1$. A sequence $(\cD_t\mid t\in T)$  is called a $\diamondsuit_T^{*}$-sequence if:
    \begin{enumerate}
        \item $\cD_t\subseteq \cP(t\down)$ and $|\cD_t|\leq\w$ for every $t\in T$.
        \item\label{condiciondeadivinanza3} For every $X\subseteq T$ we have $\{t\in T\mid X\cap t\down\not\in\cD_t\}\in NS^T$. 
    \end{enumerate}
\end{mydef}

As usual, the existence of a $\diam_T$‑sequence is denoted simply by $\diam_T$, and similarly for $\diam_T^*$. As we will see later—and in analogy with the classical case—$\diam_T^*$ implies $\diam_T$.

We begin the study of $\diamondsuit_T$ by showing that, for $\omega_1$-trees, it is at least as strong as the classical $\diamondsuit$. However, for our later work—where we separate $\diamondsuit_T$ from $\diamondsuit$—it will be useful to prove a somewhat stronger statement, one whose immediate corollary is $\diamondsuit_T \implies \diamondsuit$.

\begin{theorem}\label{improving diamT implies diam}
   Let $T$ be a nonspecial $\w_1$-tree and $S\subseteq\w_1$ such that $T\restriction S\in NS^T$, then $\diamT\implies\diam{(\w_1\setminus S)}$.
\end{theorem}
\begin{proof}
    Let \((D_t)_{t\in T}\) be a \(\diamT\)-sequence. Now, for every \(\alpha \in \omega_1\), define \(A_\alpha = \{\widehat{D_t} \mid \htop(t) = \alpha\}\). Note that \(A_\alpha \subseteq \mathcal{P}(\alpha)\) and \(|A_\alpha| \leq \omega\) for every \(\alpha \in \omega_1\).

\begin{claim}
    \((A_\alpha)_{\alpha \in \omega_1 \setminus S}\) is a \(\diam_{\omega_1 \setminus S}^{-}\)-sequence.
\end{claim}

\textit{Proof of the Claim.} Let \(X \subseteq \omega_1\), and note that \(\{t \in T \mid (T \restriction X) \cap t\down = D_t\} \notin NS^T\) since \((D_t)_{t \in T}\) is a \(\diamT\)-sequence. Now, let \(E := \{t \in T \mid (T \restriction X) \cap t\down = D_t\} \setminus (T \restriction S)\) and note that \(E \notin NS^T\). Thus, \(\widehat{E} \subseteq \{\alpha \in \omega_1 \setminus S \mid X \cap \alpha \in A_\alpha\}\). Indeed, if \(\alpha \in \widehat{E}\), then $\alpha\in\w_1\setminus S$ and there is \(t \in E\) such that \(\htop(t) = \alpha\), so \((T \restriction X) \cap t\down= D_t\), and thus \(X \cap \alpha = \widehat{D_t} \in A_\alpha\). Finally, \(\widehat{E}\) is stationary in \(\omega_1\), since otherwise \(T \restriction \widehat{E} \in NS^T\) by \ref{2 de restrictionclubnotstationary} of Theorem \ref{restrictionclubnotstationary}, but \(E \subseteq T \restriction \widehat{E}\), so \(E\) would be in \(NS^T\), which is a contradiction. \hfill\ensuremath{_{Claim}\square}
\end{proof}

\begin{corollary}\label{diamT implies diam}
    Let $T$ be a nonspecial $\w_1$-tree, then $\diamT\implies\diam$.
\end{corollary}


As mentioned earlier, Kunen proved that $\diamondsuit^-$ and $\diamondsuit$ are equivalent, which in particular allows one to deduce $\diamondsuit$ from $\diamondsuit^*$. We now establish the analogous theorem for the tree-version of these principles. 
To that end, we first prove a preliminary lemma.

\begin{lemma}\label{special bijection 2}
    Let $T$ be a well-pruned $\omega_1$-tree. Then there exists a surjection $f:T\to \omega\times T$ such that for every $t\in T$ of limit height, the restriction $f\restriction t\down$ is a bijection onto $\omega\times t\down$.
\end{lemma} 
\begin{proof}
    Let $C=\mathrm{LIM}(\omega_1)\cup \{0\}$. We will construct a sequence $\langle f_\alpha \mid \alpha \in C \rangle$ such that for every $\alpha \in C$:
    \begin{enumerate}[label=(\arabic*)]
        \item\label{cond:onto} $f_\alpha:T_{<\alpha}\to \omega\times T_{<\alpha}$ is a surjective function.
        \item\label{cond:increasing} If $\beta < \alpha$ and $\beta \in C$, then $f_\beta \subseteq f_\alpha$.
        \item\label{cond:branch} The restriction $f_\alpha\restriction t\down$ is a bijection onto $\omega\times t\down$ for every $t\in T_\alpha$.
    \end{enumerate}

    Observe that if such a sequence $\langle f_\alpha \mid \alpha \in C \rangle$ exists, then $f=\bigcup_{\alpha\in C}f_\alpha$ is a surjection from $T$ onto $\omega\times T$ satisfying the requirements of the lemma.

    We proceed by recursion. Let $f_0=\emptyset$. Now, let $\beta\in C$ and suppose that $f_\alpha$ has been defined for every $\alpha\in C\cap\beta$.

    \medskip
    \noindent \textbf{Case $\beta=\alpha+\omega$ for some $\alpha \in C$.}
    Enumerate the level $T_{\beta} = \{t^\beta_n \mid n \in \omega\}$. We construct a sequence of functions $\langle g_n^\beta \mid n \in \omega \rangle$ such that for every $n \in \omega$:
    \begin{enumerate}[label=(\roman*)]
        \item\label{subcond:dom} $\mathrm{dom}(g_n^\beta)=(\bigcup_{i\le n} t_i^\beta\down)\setminus T_{<\alpha}$ and the range of $g_n^\beta$ is contained in $\omega\times T_{<\beta}$.
        \item\label{subcond:height} $\mathrm{ht}(\pi_1(g_n^\beta(s)))\leq\mathrm{ht}(s)$ for every $s\in\mathrm{dom}(g^\beta_n)$.
        \item\label{subcond:ext} $g_n^\beta\subseteq g_{n+1}^\beta$.
        \item\label{subcond:onto_branch} The restriction $g_n^\beta \restriction (t_i^\beta\down\setminus T_{<\alpha})$ is a bijection onto $\omega\times (t_i^\beta\down\setminus T_{<\alpha})$ for all $i\leq n$.
    \end{enumerate}

    For the base step, let $g_0^\beta$ be any bijection from $t_0^\beta\down\setminus T_{<\alpha}$ onto $\omega\times (t_0^\beta\down\setminus T_{<\alpha})$ satisfying condition \ref{subcond:height}. 

    Suppose $g_n^\beta$ has been defined and let
    \[
    X = \left(t_{n+1}^\beta\down \setminus \bigcup_{i\le n}t_i^\beta\down\right) \setminus T_{<\alpha} \quad \text{and} \quad Y = \left(\bigcup_{i\le n}(t_i^\beta\down \cap t_{n+1}^\beta\down)\right) \setminus T_{<\alpha}.
    \]
    Observe that $X \cup Y = t_{n+1}^\beta\down \setminus T_{<\alpha}$ and $s < t$ for all $s \in Y$ and $t \in X$.
 
    We want to define $g_{n+1}^\beta$. If $t_{n+1}^\beta\down = t_{i}^\beta\down$ for some $i\leq n$, there is nothing to do; we simply set $g_{n+1}^\beta = g_n^\beta$ and we are done. We may therefore assume that $t_{n+1}^\beta\down \neq t_i^\beta\down$ for every $i \le n$, which implies, in particular, that $X$ is infinite and $Y$ is finite\footnote{Note that if $T$ is Hausdorff, this second case always holds.}. 

    Let $R = (\omega\times X)\cup ((\omega\times Y)\setminus g_n^\beta[Y])$ and fix a bijection $h:X \to R$ such that $\mathrm{ht}(\pi_1(h(s))) \leq \mathrm{ht}(s)$ for every $s \in X$. Note that $h[X]=R \subseteq \omega\times (t_{n+1}^\beta\down \setminus T_{<\alpha})$.    

    Define $g_{n+1}^\beta = g_n^\beta \cup h$. We claim that $g_{n+1}^\beta$ satisfies all four conditions. Conditions \ref{subcond:dom} and \ref{subcond:ext} are clear by definition. Condition \ref{subcond:height} holds for $g_{n+1}^\beta$ because it holds for both $g_n^\beta$ and $h$. It remains to verify condition \ref{subcond:onto_branch} for $i=n+1$ (for $i \le n$, it follows from the inductive hypothesis on $g_n^\beta$). Note that $g_{n+1}^\beta\restriction(t_{n+1}^\beta\down\setminus T_{<\alpha})=(g_n^\beta\restriction Y)\cup h$.
    
    \begin{claim}\label{surjectivity}
        $g_{n+1}^\beta[t_{n+1}^\beta\down \setminus T_{<\alpha}] = \omega \times (t_{n+1}^\beta\down \setminus T_{<\alpha})$.
    \end{claim}

    \noindent \textit{Proof of the claim:}
    [$\subseteq$] Let $s \in t_{n+1}^\beta\down \setminus T_{<\alpha}$. If $s \in X$, then by definition $g^\beta_{n+1}(s) = h(s) \in R\subseteq \omega \times (t_{n+1}^\beta\down \setminus T_{<\alpha})$. If $s \notin X$, then $s \in Y$, so there exists $i \le n$ such that $s \in t_i^\beta\down$. Thus $g_{n+1}^\beta(s) = g_n^\beta(s) = (m,t)$. By the inductive hypothesis (conditions \ref{subcond:dom} and \ref{subcond:height} for $g_n^\beta$), we have $\alpha \leq \mathrm{ht}(t) \leq \mathrm{ht}(s)$. By condition \ref{subcond:onto_branch}, $t \in t_i^\beta\down \setminus T_{<\alpha}$. This implies $t \leq s < t_i^\beta$. Since $s < t_{n+1}^\beta$, we have $t < t_{n+1}^\beta$. Given $\mathrm{ht}(t) \geq \alpha$, we conclude $(m,t) \in \omega \times (t_{n+1}^\beta\down \setminus T_{<\alpha})$.

    \noindent [$\supseteq$] Let $(m,t)\in \omega \times (t_{n+1}^\beta\down \setminus T_{<\alpha})$. 
    If $t\in X$, we are done as $\omega\times X\subseteq R=h[X]$. 
    So suppose $t\in Y$. 
    If $(m,t)\in g_n^\beta[Y]$, we are done as $g_n^\beta[Y]\subseteq g_{n+1}^\beta[t_{n+1}^\beta\down \setminus T_{<\alpha}]$. Otherwise, $(m,t)\in (\omega\times Y)\setminus g_n^\beta[Y]\subseteq R$, so $(m,t)\in h[X]$.  \finishclaim
    \medskip

    Claim \ref{surjectivity} implies that $g_{n+1}^\beta\restriction (t_{n+1}^\beta\down\setminus T_{<\alpha})$ is surjective onto $\omega \times (t_{n+1}^\beta\down \setminus T_{<\alpha})$. To finish, we need to prove that $g_{n+1}^\beta\restriction(t_{n+1}^\beta\down\setminus T_{<\alpha})$ is injective. 
    
    Since $g_{n+1}^\beta\restriction(t_{n+1}^\beta\down\setminus T_{<\alpha})=(g_n^\beta\restriction Y)\cup h$, and $h$ is injective by construction, we first verify that $g_n^\beta \restriction Y$ is injective. For this, note that there exists $i_0 \le n$ such that $Y \subseteq t_{i_0}^\beta\down \setminus T_{<\alpha}$ (simply take $i_0$ such that $t_{i_0}^\beta\down$ contains the maximal node of the finite chain $Y$). By the inductive hypothesis (condition \ref{subcond:onto_branch}), $g_n^\beta$ is a bijection on $t_{i_0}^\beta\down \setminus T_{<\alpha}$, and thus its restriction to the subset $Y$ is injective.

    Since both pieces ($g_n^\beta \restriction Y$ and $h$) are injective, it suffices to show that their ranges are disjoint, i.e., $h[X]\cap g_n^\beta[Y]=\emptyset$.
    
    First, observe that $g_n^\beta[Y] \cap (\omega \times X) = \emptyset$. To see this, let $s \in Y$ and $s' \in X$. Since $s < s'$, we have $\mathrm{ht}(s) < \mathrm{ht}(s')$. By condition \ref{subcond:height}, if we write $g_n^\beta(s)=(m,t)$, then $\mathrm{ht}(t) \le \mathrm{ht}(s) < \mathrm{ht}(s')$, which implies $t \neq s'$. As this holds for any $s' \in X$, we conclude that $t \notin X$, and therefore $g_n^\beta(s) \notin \omega \times X$. Second, observe that $g_n^\beta[Y]$ is disjoint from $(\omega\times Y)\setminus g_n^\beta[Y]$ simply by definition.
    Therefore, $h[X]\cap g_n^\beta[Y]=\emptyset$, as desired.

    Finally, let $f_\beta = f_\alpha \cup (\bigcup_{n \in \omega} g_n^\beta)$. Condition \ref{subcond:onto_branch} implies that $f_\beta$ satisfies condition \ref{cond:branch}. Condition \ref{cond:onto} follows from the fact that $T$ is well-pruned: for every $s \in T_{<\beta} \setminus T_{<\alpha}$, there exists $n \in \omega$ such that $s \in t_n^\beta\down$, which implies that $s \in \mathrm{dom}(f_\beta)$. Furthermore, condition \ref{subcond:onto_branch} ensures that for every such $s$ and every $m \in \omega$, the pair $(m,s)$ is in the range of $f_\beta$. Finally, condition \ref{cond:increasing} is satisfied by construction.

    \medskip
    \noindent \textbf{Case $\beta$ is a limit point of $C \cap \beta$.} In this case, we simply define $f_\beta = \bigcup_{\alpha \in C \cap \beta} f_\alpha$.
\end{proof}

\begin{mydef}
    If $T$ is a nonspecial tree of height $\w_1$. A sequence $(\cD_t\mid t\in T)$  is called a $\diamondsuit_T^{-}$-sequence if:
    \begin{enumerate}
        \item $\cD_t\subseteq \cP(t\down)$ and $|\cD_t|\leq\w$ for every $t\in T$.
        \item\label{condiciondeadivinanza2} For every $X\subseteq T$ we have $\{t\in T\mid X\cap t\down\in\cD_t\}\not\in NS^T$. 
    \end{enumerate}
\end{mydef}


\begin{theorem}\label{kunenfortrees}
  Let $T$ be a nonspecial, well-pruned $\w_1$-tree. Then $\diamT^{-}$ and $\diamT$ are equivalent.
\end{theorem}

\begin{proof}
Let $(\cD_t)_{t\in T}$ be a $\diamT^{-}$-sequence and let us fix a surjection $f:T\to\w\times T$ as in Lemma \ref{special bijection 2} and call $C=\mathrm{LIM}(\w_1)\cup\{0\}$. Now for every $t\in T$ let $\cB_t=\{f[A]\mid A\in\cD_t\}$ in case that $t\in T\restriction C$ and $\cB_t=\emptyset$ otherwise. 

\begin{claim}\label{we can guess subsests of w times the tree}
    For every $X\subseteq\w\times T$ we have that $\{t\in T\mid X\cap (\w\times t\down)\in \cB_t\}\not\in NS^T$.
\end{claim}
\textit{Proof of the Claim.} We know that $R:=\{t\in T\mid f^{-1}(X)\cap t\down\in\cD_t\}\not\in NS^T$ as $(\cD_t)_{t\in T}$ is a $\diamT^{-}$-sequence. On the other hand,  by Theorem \ref{restrictionclubnotstationary}, $T\restriction C\in (NS^T)^*$, thus $R\cap (T\restriction C)\not\in NS^T$. Now if $t\in R\cap(T\restriction C)$ we know that: $$f^{-1}(X)\cap t\down\in\cD_t\implies X\cap f[t\down]\in\cB_t\implies X\cap (\w\times t\down)\in\cB_t,$$
that is, $R\cap (T\restriction C)\subseteq\{t\in T\mid X\cap (\w\times t\down)\in \cB_t\}$, so this last is not in $NS^T$.\finishclaim

As each $\cB_t$ is countable for every $t\in T\restriction C$, let us write it as $\cB_t=\{B_t^k\mid k\in\w\}$. Clearly in this case $B_t^k\subseteq \w\times t\down$. Now, let $B_{t,n}^k=\{s\in T\mid (n,s)\in B_t^k\}$ for every $n\in\w$. 

\begin{claim}
    There is $n\in\w$ such that $(B_{t,n}^n)_{t\in T}$ is a $\diamT$-sequence.
\end{claim}
\textit{Proof of the Claim.}
Suppose that it is false, then for every $n\in\w$ there is $Q_n\subseteq T$ such that $Z_n:=\{t\in T\mid Q_n\cap t\down= B_{t,n}^n\}\in NS^T$. Now, let $Q:=\bigcup_{n\in\w}(\{n\}\times Q_n)$. For every $n\in\w$ we know that $Y_n:=\{t\in T\mid Q\cap (\w\times t\down)=B_t^n\}\in NS^T$. Indeed, to see that $Y_n\in NS^T$ it is enough to see that $Y_n\subseteq Z_n$, so:
\begin{multline*}
    t\in Y_n\implies Q\cap (\w\times t\down)=B_t^n\implies \{s\in t\down\mid (n,s)\in Q\cap(\w\times t\down)\}=B_{t,n}^n\implies\\
    \{s\in t\down\mid (s\in t\down)\wedge(s\in Q_n)\}=B_{t,n}^n\implies Q_n\cap t\down=B_{t,n}^n\implies t\in Z_n.
\end{multline*}

Now, since the ideal $NS^T$ is $\w_1$-complete (part \ref{ideal of non-stationary is complete} of Theorem \ref{teo que junta cosas sobre NS^T}), we know that the set $\bigcup_{n\in\w}Y_n=\{t\in T\mid \exists n\in\w(B\cap (\w\times t\down)=B_t^n)\}\in NS^T$, which is a contradiction to Claim \ref{we can guess subsests of w times the tree}.
\hfill\ensuremath{_{Claim}\square}
\end{proof}

As we said before, a direct consequence of Theorem~\ref{kunenfortrees} is that $\diam_T^* \Longrightarrow \diam_T$ for every nonspecial well-pruned $\omega_1$‑tree $T$.

What we consider the most natural next step is to explore how $\diamondsuit_S$ and $\diamondsuit_T$ compare when $S$ is a nonspecial subtree of $T$. For this, it is convenient to recall the natural ordering on arbitrary trees:
\begin{notation}
If $S$ and $T$ are trees, we write $S\le T$ to mean that there exists a strictly increasing map $g\colon S\to T$, i.e., such that $g(s)<_T g(t)$ whenever $s<_S t$.
\end{notation}

Observe that if $g \colon S \to T$ is strictly increasing, then the map $f(s) = g(s) \upharpoonright \operatorname{ht}(s)$ is also strictly increasing. As defined in \cite{stevo-lipschitzmaps}, a strictly increasing and level-preserving\footnote{A map \(f\colon S\to T\) is \emph{level-preserving} if \(\operatorname{ht}(s)=\operatorname{ht}(f(s))\) for every \(s\in S\).} map \(f\colon S\to T\) is called a \emph{Lipschitz map}\footnote{Todorcevic's notion is more general as it allows partial maps, but when the domain is downward closed (in particular, when it is all of \(S\)) it is equivalent to being strictly increasing and level-preserving.}. Hence, whenever \(S\le T\) and \(f\colon S\to T\) witnesses this relation, we will assume that \(f\) is Lipschitz.

However, a Lipschitz map \(f\colon S\to T\) need not be injective: it may send distinct nodes of the same level in \(S\) to a single node in \(T\). For instance, for any tree \(T\) of height $\w_1$ we have \(T\leq\omega_1\) via the map sending each node to its height. Conversely, \(T\) has a cofinal branch precisely when \(\omega_1\leq T\).

Note that if \(f\colon S\to T\) is Lipschitz and \(A\subseteq T\) is an antichain, then \(f^{-1}(A)\) is also an antichain. Indeed, if \(s<_S t\) both lie in \(f^{-1}(A)\), then \(f(s)<_T f(t)\) and both are in \(A\), contradicting that \(A\) is an antichain. As a consequence, the preimage of any special subset \(B\subseteq T\) is special; in particular we have the following fact.
\begin{fact}\label{Lipschitz map of a nonspecial is nonspecial}
If \(S\le T\) and \(S\) is nonspecial, then \(T\) is nonspecial.
\end{fact}

\begin{notation}
Let $s\in S$ and $\alpha\le\htop(s)$.  We denote by $s\restriction\alpha$ the unique predecessor of $s$ with height $\alpha$, that is, the unique $u\in S$ with $u<_S s$ and $\htop(u)=\alpha$.
\end{notation}

\begin{fact}\label{Lipschitz map commutes with restriction}
Let $f\colon S\to T$ be a Lipschitz map and let $s\in S$.  Then for every $\alpha\le\htop(s)$,
\[
f\bigl(s\restriction\alpha\bigr)
\;=\;
f(s)\restriction\alpha.
\]
\end{fact}

\begin{proof}
Set $v = s\restriction\alpha$.  Note that $f(v)$ satisfies:
\begin{itemize}
  \item $f(v)<_T f(s)$ since $v<_S s$, and
  \item $\htop(f(v)) = \alpha = \htop(f(s)\restriction\alpha)$.
\end{itemize}
Since $f(s)\restriction\alpha$ is the unique node in $T$ with these properties, we necessarily have $f(v) = f(s)\restriction\alpha$.
\end{proof}

In order to state the following result, we first need to introduce the \Katetov order.

\begin{mydef} \cite{katetovproduct} Let $\cI$ and $\cJ$ be ideals on $X$ and $Y$ respectively and $f:Y\rightarrow X$. \begin{enumerate} \item $f$ is a \emph{\Katetov}function from $(Y,\cJ)$ to $(X,\cI)$ if $f^{-1}(A)\in\cJ$ for all $A\in\cI$. \item $\cI\leq_K\cJ$ (\emph{$\cI$ is \Katetov below $\cJ$}) if there exists a \Katetov function from $(Y,\cJ)$ to $(X,\cI)$. \end{enumerate} \end{mydef}

The interested reader is referred to \cite{michaelcombinatoricoffilters,michaelkatetov} for further details on the Katětov order.

\begin{theorem}\label{Lipschitz map and katetov}
Let $S$ and $T$ be trees of height $\omega_1$ with $S\leq T$.  Then:
\begin{enumerate}[label=(\roman*)]
  \item\label{a de gfgfgf} $NS^T\le_K NS^S$, and
  \item\label{b de gfgfgf} If $S$ is a
  nonspecial $\w_1$-tree\footnote{By Fact \ref{Lipschitz map of a nonspecial is nonspecial}, this implies that $T$ is also nonspecial, in particular it makes sense to consider $\diamT$.}, then $\diam_S$ implies $\diam_T$.
\end{enumerate}   
\end{theorem}

\begin{proof}
   Let $f\colon S\to T$ be a Lipschitz map. To prove \ref{a de gfgfgf}, we show that $f$ is a Katětov function from $(T,NS^T)$ to $(S,NS^S)$, i.e., whenever $X\subseteq T$ satisfies $X\in NS^T$, then $f^{-1}(X)\in NS^S$. Fix such an $X\in NS^T$. By definition,
$$
X=\bigcup_{t\in T}\bigl(A_t\cap t\up\bigr),
$$
where each $A_t$ is special, and hence each $f^{-1}(A_t)$ is also special.

For each $s\in S$, set $B_s = f^{-1}\bigl(A_{f(s)}\bigr)$, which is special. We will show
$$
f^{-1}(X)\;\subseteq\;\bigcup_{s\in S}\bigl(B_s\cap s\up\bigr),
$$
and since $\bigcup_{s\in S}(B_s\cap s\up)\in NS^S$, this will prove $f^{-1}(X)\in NS^S$.

Take any $u\in f^{-1}(X)$. Then $f(u)\in X=\bigcup_{t\in T}(A_t\cap t\up)$, so pick $t\in T$ with $f(u)\in A_t\cap t\up$.  Let $\alpha=\htop(t)$ and $s = u\restriction\alpha$.  By Proposition~\ref{Lipschitz map commutes with restriction}, $f(s)=t$.  Now
\[
u\in B_s \;\Longleftrightarrow\; u\in f^{-1}(A_{f(s)})\;\Longleftrightarrow\; u\in f^{-1}(A_t)
\;\Longleftrightarrow\; f(u)\in A_t,
\]
where the final condition holds by choice of $t$. Clearly $u\in s\up$, and hence $u\in B_s\cap s\up$, as required.

To prove \ref{b de gfgfgf} assume $S$ is an $\omega_1$-tree and let $(D_s)_{s\in S}$ be a $\diam_S$-sequence.  For each $t\in T$, set
\[
\mathcal D_t = \{\,f[D_s]\mid s\in S,\ f(s)=t\}.
\]
Note that each $\mathcal D_t$ is a subfamily of $\mathcal P(t\down)$ and is countable, since $f^{-1}(\{t\})\subseteq S_{\htop(t)}$ and $S_{\htop(t)}$ is countable because $S$ is an $\omega_1$-tree.

We will show that $(\mathcal D_t)_{t\in T}$ is a $\diam^-_T$-sequence.  Fix $X\subseteq T$.

\begin{claim}\label{algun claim z}
For any $s\in S$, we have $f\bigl(s\down\cap f^{-1}(X)\bigr)\;=\;f(s)\down\cap X$.
\end{claim}

\textit{Proof of Claim:}
“$\subseteq$”: If $u\in s\down\cap f^{-1}(X)$ then $u<_S s$ and $f(u)\in X$.  Since $f$ is a Lipschitz map, $f(u)<_T f(s)$, so $f(u)\in f(s)\down\cap X$.  

“$\supseteq$”: Let $w\in f(s)\down\cap X$ and call $v=s\restriction\htop(w)$. By Fact \ref{Lipschitz map commutes with restriction}, we have that $f(v)=f(s)\restriction\htop(w)=w$.  Since $w\in X$, it follows that $v\in f^{-1}(X)$, and hence $v\in s\down\cap f^{-1}(X)$, proving $w\in f\bigl(s\down\cap f^{-1}(X)\bigr)\,$.
\finishclaim

Let 
\[
Y \;=\;\{\,t\in T \mid X\cap t\down\in\mathcal D_t\}.
\]
We must show $Y\notin NS^T$.  Suppose otherwise.  Then by \ref{a de gfgfgf}, $f^{-1}(Y)\in NS^S$.  On the other hand, since $(D_s)$ is a $\diam_S$-sequence,
\[
Z \;=\;\{\,s\in S \mid f^{-1}(X)\cap s\down = D_s\}\not\in NS^S.
\]
But Claim~\ref{algun claim z} implies that 
\[
s\in Z
\;\Longrightarrow\;
f(s)\in Y,
\]
so $Z\subseteq f^{-1}(Y)$, contradicting the nonstationarity of $f^{-1}(Y)$.  Hence $Y\notin NS^T$, and $(\mathcal D_t)_{t\in T}$ is indeed a $\diam^-_T$-sequence. Finally, by Theorem \ref{kunenfortrees}, we conclude that $\diam_T$ holds.
\end{proof}

Note that part (ii) of Theorem~\ref{Lipschitz map and katetov} provides an alternative proof of Corollary~\ref{diamT implies diam}, since if $T$ is a nonspecial $\omega_1$-tree, then $T \leq \omega_1$. Another consequence of Theorem~\ref{Lipschitz map and katetov} is the following.

\begin{corollary}\label{diamond implies diamond for trees with cofinal branch}\label{diamond goes up}
Let $S$ and $T$ be nonspecial trees of height~$\omega_1$ such that $S$ is a subtree of $T$. Then $\diamondsuit_S$ implies $\diamondsuit_T$.
\end{corollary}

\begin{proof}
The inclusion $S\hookrightarrow T$ is a Lipschitz map, so part (ii) of Theorem~\ref{Lipschitz map and katetov} gives $\diamondsuit_S\implies\diamondsuit_T$.  
\end{proof}

A particularly interesting case is when $S$ is a cofinal branch of $T$. Note that in this scenario, we have $\w_1 \leq S \leq \w_1$, which implies that $S$ is a nonspecial $\w_1$-tree and, by part (ii) of Theorem \ref{Lipschitz map and katetov}, $\diam_S$ and $\diam$ are equivalent. In particular, by Corollary \ref{diamond goes up}, we have that $\diam$ implies $\diamT$. This, combined with Corollary \ref{diamT implies diam}, yields the following result.

\begin{corollary}\label{coro 3.21}
For a tree $T$ of height $\w_1$ with a cofinal branch, $\diam$ implies $\diamT$. Additionally, if $T$ is an $\w_1$-tree, then $\diam$ and $\diamT$ are equivalent.
\end{corollary}

An important nonspecial $\w_1$-tree is
\[
T(\emptyset)=\{s\in 2^{<\w_1}\mid |\operatorname{supp}(s)|<\w\},
\]
which possesses several interesting properties. In particular, it is shown in \cite[Theorem 4.1 and Corollary 4.4]{stevotreesandlinearly} that $T(\emptyset)$ has no Aronszajn subtrees and is \emph{almost-Suslin} (Definition \ref{def de almost Suslin}).  A natural question is what can be said about $\diam_{T(\emptyset)}$. Corollary \ref{coro 3.21} answers it: since $T(\emptyset)$ has a cofinal branch, $\diam_{T(\emptyset)}$ is equivalent to the usual $\diam$; indeed, as we will show in Theorem \ref{diamT sii diam for almost-suslin}, being almost-Suslin is a sufficient condition for this equivalence.

We now turn to a different topic. One of the earliest results, if not the very first, when studying $\diamondsuit$ is that it implies the Continuum Hypothesis ($\mathsf{CH}$). Therefore, when studying $\diamondsuit_T$, it is natural to ask what cardinal arithmetic constraints $\diamondsuit_T$ imposes. Of course, by Corollary \ref{diamT implies diam}, we know that $\diamondsuit_T$ implies $\mathsf{CH}$ whenever $T$ is an $\omega_1$-tree. However, it is interesting to explore whether $\diamondsuit_T$ generally provides a cardinal inequality associated with the levels of the tree $T$.

Recall that the standard proof that $\diamondsuit$ implies $\mathsf{CH}$ shows that if $(D_\alpha)_{\alpha < \omega_1}$ is a $\diamondsuit$-sequence, then for every subset $X \subseteq \omega$, there exists some $\alpha \in \omega_1 \setminus \omega$ with $X = D_\alpha$. Consequently,
\[
|\mathcal{P}(\omega)| \leq |\omega_1 \setminus \omega| = \omega_1,
\]
so $\mathsf{CH}$ follows.

In other words, when we think of $\omega_1$ as a tree $T$, the above can be reformulated as $\diamondsuit_T$ implying that $|\mathcal{P}(T_{< \omega})| \leq |T_{\geq \omega}|$. It is natural to ask whether this inequality holds for any tree $T$ of height $\omega_1$ whenever $\diamondsuit_T$ holds. However, this is not true in general. For example, consider a tree $S$ of height $\omega$ with a cofinal branch $\{r_n \mid n \in \omega\}$ and such that $|S| = 2^{\omega_1}$. Now, let $R = \{r_\alpha \mid \alpha \in \omega_1 \setminus \omega\}$ be a set disjoint from $S$, where $r_\alpha \neq r_\beta$ if $\alpha \neq \beta$. Define $T = S \cup R$ as the tree where $r_\alpha \down = \{r_\beta \mid \beta < \alpha\}$ and that has $S$ a subtree. Suppose $\diamondsuit$ holds; since $T$ has a cofinal branch, by Corollary \ref{coro 3.21}, $\diamondsuit_T$ also holds. However, note that 
\[
|T_{\geq \omega}| = \omega_1 < |T| = 2^{\omega_1} < |\mathcal{P}(T_{< \omega})| = |\mathcal{P}(S)| = |\mathcal{P}(2^{\omega_1})| = 2^{2^{\omega_1}}.
\]
This yields the following result:

\begin{proposition}\label{prop:cardinality_not_like_w1}
    It is consistent that there exists a nonspecial tree $T$ of height $\omega_1$ such that $\diamondsuit_T$ holds and $|T_{\geq \omega}| < |T| < |\mathcal{P}(T_{< \omega})|$.
\end{proposition}

Although $|\mathcal{P}(T_{< \omega})| \leq |T_{\geq \omega}|$ does not necessarily hold whenever $\diamondsuit_T$ holds, we can still prove the following:

\begin{proposition}\label{prop:cardinality_like_w1}
    Let $T$ be a nonspecial tree of height $\omega_1$. Then $\diamondsuit_T$ implies $2^\omega \leq |T_{\geq \omega}|$.
\end{proposition}
\begin{proof}
    Let $(D_t)_{t \in T}$ be a $\diamondsuit_T$-sequence. By Theorem \ref{restrictionclubnotstationary}, $T_{\geq \omega} = T \restriction (\omega_1 \setminus \omega) \in (NS^T)^*$. Now fix $A \subseteq \omega$. Since $(D_t)_{t \in T}$ is a $\diamondsuit_T$-sequence, there exists $t_A \in T_{\geq \omega}$ such that $(T \restriction A) \cap t_A \down = D_{t_A}$, i.e., $\{s < t_A \mid \htop(s) \in A\} = D_{t_A}$. It is straightforward to verify that $t_A \neq t_B$ if $A \neq B$. Therefore, $2^\omega \leq |T_{\geq \omega}|$.
\end{proof}

In particular, as a consequence of Proposition \ref{prop:cardinality_like_w1}, we have an alternative proof of the fact that if $T$ is a nonspecial $\omega_1$-tree, then $\diamondsuit_T$ implies $\mathsf{CH}$, since in that case $|T_{\geq \omega}| = \omega_1$.

Given the flexibility of the definition of $\diam_T$, which allows us to consider any tree of height $\omega_1$, a natural question arises: does there exist a tree $T$ such that $\diam_T$ is provable in $\mathsf{ZFC}$? We know that for $T=\omega_1$ the answer is negative, since $\diamondsuit$ is independent of $\mathsf{ZFC}$.

However, as we will see, $2^{<\omega_1}$ is a tree whose diamond principle holds in $\mathsf{ZFC}$. To prepare for that argument, it is convenient to fix some notation.

\begin{notation}
If $T$ is a tree and $t\in T$, then $\mathrm{Succ}_T(t)$ denotes the set
$$ \{\,s\in T \mid t<_T s \text{ and } \htop(s)=\htop(t)+1\,\},$$
that is, the collection of immediate successors of $t$ in $T$.
\end{notation}

\begin{mydef}
Let $T$ be a tree.
\begin{enumerate}
  \item A family $\mathcal Q=\{Q_t\mid t\in T\}$ is a \emph{successor-partition of $T$} if for every $t\in T$ the set \mbox{$Q_t=\{Q_t(0),Q_t(1)\}$} is a partition of $\mathrm{Succ}_T(t)$; i.e., $Q_t(0)\cup Q_t(1)=\mathrm{Succ}_T(t)$ and $Q_t(0)\cap Q_t(1)=\emptyset$.
  \item Given a successor-partition $\mathcal Q$ of $T$ and a function $f\colon T\to 2$, the subtree of $T$ determined by $\mathcal Q$ and $f$, denoted $T(\mathcal Q,f)$, is the unique subtree of $T$ satisfying:
    \begin{enumerate}
      \item If $t\in T(\mathcal Q,f)$, then $\mathrm{Succ}_{T(\mathcal Q,f)}(t)=Q_t({f(t)})$.
      \item If $t\in T$, $\htop(t)$ is a limit ordinal and $t\restriction\beta\in T(\mathcal Q,f)$ for every $\beta<\htop(t)$, then $t\in T(\mathcal Q,f)$.
    \end{enumerate}
  \item A successor-partition $\mathcal Q$ of $T$ is \emph{nice} if $T(\mathcal Q,f)$ is nonspecial for every $f\colon T\to 2$.
\end{enumerate}
\end{mydef}

\begin{theorem}\label{nice succesor implies diamond}
Let $T$ be a tree of height $\w_1$. If $T$ admits a nice successor-partition, then $\diam_T$ holds.
\end{theorem}

\begin{proof}
For each $t\in T$ define
$D_t=\{\,t\restriction\alpha \mid \alpha<\htop(t)\text{ and }t\restriction(\alpha+1)\in Q_{\,t\restriction\alpha}(1)\,\}$.
For any $X\subseteq T$ let $\chi_X\colon T\to 2$ be the characteristic function of $X$.
Since $\mathcal Q$ is nice, the subtree $S:=T(\mathcal Q,\chi_X)$ is nonspecial.

\begin{claim}
$S\subseteq\{\,t\in T \mid X\cap t\downarrow = D_t\,\}$.
\end{claim}

\textit{Proof of the claim.}
Fix $t\in S$ and $\alpha<\htop(t)$. We show that
$t\restriction\alpha\in D_t$ if and only if $t\restriction\alpha\in X$.

By definition,
\begin{equation}\label{ccccc}
 t\restriction\alpha\in D_t \iff t\restriction(\alpha+1)\in Q_{\,t\restriction\alpha}(1)   
\end{equation}
On the other hand, because $S$ is a subtree of $T$, $t\in S$ and $t\restriction(\alpha+1)\leq_T t$, we also have
$t\restriction(\alpha+1)\in S$, hence
\begin{equation}\label{eeeee}
 t\restriction(\alpha+1)\in Q_{\,t\restriction\alpha}({\chi_X(t\restriction\alpha)})   
\end{equation}
Therefore, by (\ref{ccccc}) and (\ref{eeeee}), we have 
$$t\restriction(\alpha+1)\in Q_{\,t\restriction\alpha}(1)\iff\chi_X(t\restriction\alpha)=1.$$ This yields
$t\restriction\alpha\in D_t \iff \chi_X(t\restriction\alpha)=1$, i.e.,\ 
$t\restriction\alpha\in D_t \iff t\restriction\alpha\in X$.
Consequently $X\cap t\downarrow = D_t$, as required.\finishclaim

Since $S$ is nonspecial and it is a subtree of $T$, Corollary~\ref{for subtrees special and nonstationary is the same} yields $S \notin NS^T$. Hence
\[
\{\,t\in T \mid X\cap t\downarrow = D_t\,\} \notin NS^T.
\]
As $X$ was arbitrary, the sequence $(D_t)_{t\in T}$ is a witness of $\diam_T$.
\end{proof}

\begin{lemma}\label{2w1 admists a nice succesor}
    $2^{<\omega_1}$ admits a nice successor-partition.
\end{lemma}
\begin{proof}
    Let $T=2^{<\omega_1}$. For each $t\in T$ define
    $Q_t(0)=\{t\conc 0\}$ and $Q_t(1)=\{t\conc 1\}$ and set \mbox{$\cQ=\{Q_t\mid t\in T\}$}. Fix $f\colon T\to 2$ and consider the subtree $T(\cQ,f)$.  Define a sequence $(s_\alpha)_{\alpha<\omega_1}$ by
    $s_0=\emptyset$, 
    $s_{\alpha+1}=s_\alpha\conc f(s_\alpha)$, and for limit $\alpha$ set $s_\alpha=\bigcup_{\beta<\alpha}s_\beta$.
    Then $T(\cQ,f)=\{s_\alpha\mid \alpha<\omega_1\}$, so $T(\cQ,f)$ is a cofinal branch of $T$, in particular $T(\cQ,f)$ is nonspecial. Therefore $\cQ$ is nice.
\end{proof}

As a consequence of Theorem \ref{nice succesor implies diamond} and Lemma \ref{2w1 admists a nice succesor}, we have:

\begin{theorem}\label{diam for 2^<w1 is true in ZFC}
    $\diam_{2^{<\omega_1}}$ holds in $\mathsf{ZFC}$.
\end{theorem}

Two canonical examples of nonspecial $\w_1$-trees without cofinal branches are the tree $\sigma\mathbb{Q}$ (first studied by Kurepa \cite{kurepa1936ensembles,kurepathesis}) and the ``shooting-a-club'' tree $T(S)$ (where $S$ is a bistationary subset of $\w_1$). The latter was first introduced by Baumgartner, Harrington, and Kleinberg in \cite{baumgartneraddingaclub} as a forcing notion and was later studied by Todorcevic as a tree in \cite{stevostationary}, where several interesting facts were established; for example, $T(S)\times T(S')$ is special if and only if $S\cap S'\in NS_{\w_1}$, which settles Hedetniemi's Conjecture at that level.

The tree $\sigma\mathbb{Q}$ is the set $\{\,t \mid \exists\alpha<\omega_1\ (t\colon\alpha\to\mathbb{Q}\text{ increasing and }\mathrm{ran}(t)\text{ is bounded})\,\}$ ordered by $\subseteq$. 
A natural successor-partition of $\sigma\mathbb{Q}$ is obtained by fixing a partition $\mathbb{Q}=\mathbb{Q}_0\cup\mathbb{Q}_1$ in which both $\mathbb{Q}_0$ and $\mathbb{Q}_1$ are dense in $\mathbb{Q}$. For each $t\in\sigma\mathbb{Q}$ and $i\in\{0,1\}$ set $
Q_t(i)=\{\,s\in\mathrm{Succ}_{\sigma\mathbb{Q}}(t)\mid s(\mathrm{ht}(t))\in\mathbb{Q}_i\,\}$. It is natural to ask whether the successor-partition $\mathcal{Q}$ is nice; that is, whether for every $f\colon\sigma\mathbb{Q}\to 2$, the subtree $\sigma\mathbb{Q}(\mathcal{Q},f)$ is nonspecial. The standard proof that $\sigma\mathbb{Q}$ is nonspecial, as shown in \cite[Corollary 9.9]{stevotreesandlinearly}, does not seem to immediately adapt to prove $\sigma\mathbb{Q}(\mathcal{Q},f)$ is nonspecial.

\begin{question}
    Is the successor-partition $\mathcal Q$ described above nice? Or more generally: does $\diam_{\sigma\mathbb{Q}}$ hold?
\end{question}

The ``shooting-a-club'' tree is defined as follows.  Fix a bistationary set $S\subseteq\omega_1$.  Let $T(S)$ be the collection of bounded closed subsets $p\subseteq S$ of ordinals, ordered by end-extension: $q\le p$ iff $q=p\cap\alpha$ for some $\alpha<\omega_1$.  It is straightforward to verify that $T(S)$ is well-pruned, has no cofinal branch, and is nonspecial.  The last assertion follows from the fact that any special well-pruned tree of height $\omega_1$ collapses $\omega_1$, whereas $T(S)$ does not collapse $\omega_1$; in fact $T(S)$ is $\omega$-distributive (see \cite[Theorem B]{baumgartneraddingaclub} or \cite[Lemma 23.9, p.\ 445]{jechbook}).  As with $\sigma\mathbb{Q}$, it is natural to ask:



\begin{question}
    Given a bistationary $S\subseteq\omega_1$, does $\diam_{T(S)}$ hold?
\end{question}

Finally, we conclude the section by showing that $\diamondsuit_T$ can be expressed within a more familiar framework, a connection that Assaf Rinot brought to our attention, and for which we express our gratitude.

\begin{mydef}\label{diam de ideales}
    Let $\cI$ be an ideal over the powerset of some set $Z$, i.e., $\cI\subseteq\cP(\cP(Z))$. We say that $\diam(\cI)$ holds if there exists a function $g:\cP(Z)\to \cP(Z)$ such that for every $X\subseteq Z$, the collection $\{Y\subseteq Z\mid X\cap Y=g(Y)\}$ is in $\cI^+$. 
\end{mydef}

Now suppose that $\diam(\cI)$ holds and let $g:\cP(Z)\to\cP(Z)$ be a witness to this. Now consider $f:\cP(Z)\to\cP(Z)$ defined by $f(Y)=g(Y)\cap Y$. Then $f$ is also a witness of $\diam(\cI)$. Indeed, it suffices to note that for every $X\subseteq Z$, we have:
    \[
    \{Y\subseteq Z\mid X\cap Y=g(Y)\}=\{Y\subseteq Z\mid X\cap Y=g(Y)\cap Y\}=\{Y\subseteq Z\mid X\cap Y=f(Y)\}.
    \]
Thus, without loss of generality, we can assume that $g(Y)\subseteq Y$ for all $Y\subseteq Z$.

The language of Definition \ref{diam de ideales} is usually used to define the diamond principle on $P_\kappa(\lambda)$ (see, for example, \cite{SHIOYA}), but it also captures $\diamT$, as we will see.

Let $T$ be a tree of height $\omega_1$, and let $\cJ_T\subseteq \cP(\cP(T))$ be the ideal defined as follows:
\begin{enumerate}
    \item $\{Y\subseteq T\mid\forall t\in T(t\down\neq Y)\}\in\cJ_T$, i.e., $\{t\down\mid t\in T\}\in(\cJ_T)^*$.
    \item If $A\subseteq\{t\down\mid t\in T\}$, then $A\in\cJ_T$ if and only if $\{t\in T\mid t\down\in A\}\in NS^T$.
\end{enumerate}

Essentially, $\cJ_T$ is designed to mimic the ideal $NS^T$ but over $\cP(T)$.

\begin{proposition}\label{diamT implies diam(I) for some I}
    $\diamT$ is equivalent to $\diam(\cJ_T)$.
\end{proposition}
\begin{proof}
    Suppose that $(D_t)_{t\in T}$ is a $\diamT$-sequence, and let $g:\cP(T)\to\cP(T)$ be defined by $g(t\down)=D_t$ for all $t\in T$, and constant $\emptyset$ for the rest of $\cP(T)$. It is clear that for every $X\subseteq T$, we have:
    \[
    \{Y\subseteq T\mid X\cap Y=g(Y)\}\supseteq\{t\down\mid X\cap t\down=D_t\}
    \]
    and the latter set is in $\cJ_T^+$ since $\{t\mid X\cap t\down=D_t\}\notin NS^T$, as $(D_t)_{t\in T}$ is a $\diamT$-sequence.

    Conversely, suppose that $g$ is a witness of $\diam(\cJ_T)$ such that $g(Y)\subseteq Y$  for every $Y\subseteq T$. Now, for every $t\in T$, let $D_t=g(t\down)$. Note that since $\{Y\subseteq T\mid\forall t\in T(t\down\neq Y)\}\in\cJ_T$, for every $X$ we have $\{t\down\mid X\cap t\down=D_t\}\notin\cJ_T$. Therefore, $\{t\mid X\cap t\down=D_t\}\notin NS^T$, which proves that $(D_t)_{t\in T}$ is a $\diamT$-sequence.
\end{proof}


\section{\texorpdfstring{The consistency of $\diam_T$ and the case for almost-Suslin trees}{The consistency of Diamond of a tree and the case for almost-Suslin trees}}

Although we already know some results regarding
$\diamT$, we still lack a crucial piece of information: its consistency. The first thing we will do in this section is to establish that consistency for some class of trees.

\begin{theorem}\label{implication of diamons *}
    Let $T$ be a nonspecial $\w_1$-tree. Then $\diam^*$ implies $\diamT^*$.
\end{theorem}
\begin{proof}
    Let $(\cA_\alpha\mid\alpha\in\w_1)$ be a $\diam^*$ -sequence and enumerate $\cA_\alpha=\{A_\alpha^n\mid n\in\w\}$ for every $\alpha\in\w_1$. Now assume without loss of generality that $(T,<_T)=(\w_1,<_T)$. We know that the set $$C=\{\alpha\in\w_1\mid T_{<\alpha}=\alpha\}$$
    is club, consequently, by condition \ref{4 de restrictionclubnotstationary} of Lemma \ref{restrictionclubnotstationary}, $T\restriction C\not\in NS^T$. Now for every $t\in T\restriction C$ let $\cD_t=\{A_\alpha^n\cap t\down\mid n\in\w\}$ and if $t\not\in T\restriction C$ let $\cD_t=\emptyset$. 

    \begin{claim}
       $(\cD_t\mid t\in T)$ is a $\diamT^*$-sequence. 
    \end{claim}
    \textit{Proof of the Claim.}
    Let $X\subseteq T$. As $T=\w_1$, $D=\{\alpha\in\w_1\mid X\cap\alpha\in\cA_\alpha\}$ contains a club, hence so does $E=D\cap C$. Let any $t\in T\restriction E$ and let $\alpha=\htop(t)$. Now we know that $X\cap T_{<\alpha}=X\cap\alpha\in\cA_\alpha$, so $X\cap\alpha=A_\alpha^n$ for some $n\in\w$, 
thus: $$X\cap t\down=(X\cap\alpha)\cap t\down=A_\alpha^n\cap t\down\in\cD_t.$$
    This proves that $\{t\in T\mid X\cap t\down\not\in\cD_T\}\subseteq T\restriction (\w_1\setminus E)=T\setminus (T\restriction E)\in NS^T$ as $(T\restriction E)\in (NS^T)^*$ by condition \ref{5 de restrictionclubnotstationary} of Lemma \ref{restrictionclubnotstationary}. \hfill\ensuremath{_{Claim}\square}
\end{proof}

It is natural to ask whether the proof of Theorem \ref{implication of diamons *} can be adapted to show that $\diamondsuit$ implies $\diamondsuit_T$. However, note that this seems unlikely because in the proof we use the fact that if $D \subseteq \omega_1$ (the set where $X$ is guessed) is a club, then two things happen:
\begin{enumerate}[label=(\roman*)]
    \item\label{1 of why one club proof is not a stationary proof} $E=D \cap C$ is club, and
    \item\label{2 of why one club proof is not a stationary proof} $T \restriction E \in \left(NS^T\right)^*$.
\end{enumerate}
Note that if $D$, the set where $X$ is guessed, is merely stationary (not necessarily club), although condition \ref{1 of why one club proof is not a stationary proof} has an adequate analogue (as $E = D \cap C$ is stationary), the analogue of condition \ref{2 of why one club proof is not a stationary proof} fails, since we cannot guarantee that $T \restriction E \notin NS^T$. As we will see later, this is actually the only obstruction to deducing $\diamondsuit_T$ from $\diamondsuit$ (Lemma \ref{suficient condition of T for diam implying diamT}). However, since in principle not every well-pruned $\omega_1$-tree satisfies this property (as Theorem \ref{important theorem of Shelah} witness), in order to prove that $\diamondsuit_T$ is consistent for every well-pruned nonspecial $\omega_1$-tree, we must go further and obtain $\diamondsuit_T$ not as a consequence of $\diamondsuit$ but from $\diamondsuit^*$:

\begin{corollary}\label{diamT is consistent}
Let $T$ be a well-pruned nonspecial $\omega_1$-tree. Then $\diam^*$ implies $\diam_T$. 
\end{corollary}
\begin{proof}
$\diamondsuit^*$ implies $\diamondsuit_T^*$ by Theorem~\ref{implication of diamons *}, while $\diamondsuit_T^*$ implies $\diamondsuit_T$ by Theorem~\ref{kunenfortrees}. 
\end{proof}

In particular:

\begin{corollary}\label{L modela todos los diamantes}
    $L \models \forall T \text{ nonspecial well-pruned } \omega_1\text{-tree}, \diamT \text{ holds}$.
\end{corollary}

It would be interesting to know whether a direct proof of Corollary \ref{L modela todos los diamantes} can be given using the fine structure of $L$ or the Condensation Lemma.

We now aim to prove that for every well-pruned, nonspecial Aronszajn tree $T$, $\diam_T$ is consistent. By Corollary~\ref{diamT is consistent}, this reduces to showing that one can force $\diam^*$ while preserving the fact that $T$ remains a nonspecial Aronszajn tree (being well-pruned is preserved trivially). To this end we need the following lemmas, due to Silver and Todorcevic respectively, which will also be used in subsequent sections.

\begin{lemma}{\cite[p.~387]{silverkurepaconjecture}}\label{being aronszajn is preserved by sigma closed}
Let $T$ be a well-pruned Aronszajn tree and let $\mathbb{P}$ be a $\sigma$-closed forcing. Then
$ \mathbb{P}\forces ``T \text{ is a well-pruned Aronszajn tree}."$
\end{lemma}

\begin{lemma}\cite[Lemma~12]{stevoonaconjectureofrado}\label{sigma closed does not specialize w1-trees}
Let $T$ be a nonspecial $\omega_1$-tree and let $\mathbb{P}$ be a $\sigma$-closed forcing. Then $\mathbb{P}$ does not specialize $T$.
\end{lemma}






We can now prove the desired result:

\begin{theorem}
    Let $T$ be a well-pruned, nonspecial Aronszajn tree. Then in some cofinality-preserving forcing extension, $T$ is a well-pruned nonspecial Aronszajn tree and $\diamT$ holds.
\end{theorem}
\begin{proof}
As previously stated, it suffices to show that we can force $\diam^*$ with a cofinality-preserving forcing that preserves $T$ as a well-pruned, nonspecial Aronszajn tree. By \cite[Exercise H20, p.\ 249]{kunenbook}, under $\mathsf{CH}$, there exists a $\sigma$-closed forcing $\mathbb{Q}$ such that $\mathbb{Q} \forces ``\diam^*"$ (this forcing plays a key role in Section~\ref{section diamond star es mas fuerte que diamT} and additionally forces a stronger diamond principle, but this property is sufficient for now). Furthermore, it is well-known that $\mathbb{P} = (2^{<\omega_1}, \supseteq)$ is $\sigma$-closed and forces $\mathsf{CH}$. Thus, $\mathbb{P} \star\dot{\mathbb{Q}}$ is a $\sigma$-closed forcing that forces $\diam^*$. By Lemmas~\ref{being aronszajn is preserved by sigma closed} and \ref{sigma closed does not specialize w1-trees}, $T$ remains a well-pruned, nonspecial Aronszajn tree in $V[G]$. Hence, $V[G] \models \diam_T$.
\end{proof}


It is natural to think that for trees that are not too wide, $\diamT$ should be very similar to $\diam$. As we mentioned earlier, the issue with trying to deduce $\diamT$ from $\diam$ is that we do not know that $T \restriction X \notin NS^T$ when $X \subseteq \omega_1$ is stationary. As we will see, that this is actually the only obstruction (Lemma \ref{suficient condition of T for diam implying diamT}), and that almost-Suslin trees (Definition \ref{def de almost Suslin}) do not have this obstruction (Lemma \ref{almost-Suslin satisfies the suficient condition}).

\begin{lemma}\label{suficient condition of T for diam implying diamT}
    Let $T$ be a nonspecial $\w_1$-tree such that for every $X\subseteq\w_1$, $T\restriction X\in NS^T$ if and only if $X$ is nonstationary. Then $\diam$ and $\diamT$ are equivalent.
\end{lemma}
\begin{proof}
   By Corollary \ref{diamT implies diam}, we know that $\diamT$ implies $\diam$. Now suppose that $(A_\alpha\mid\alpha\in\w_1)$ is a $\diam$-sequence and assume without loss of generality that $(T,<_T)=(\w_1,<_T)$. We know that the set $$C=\{\alpha\in\w_1\mid T_{<\alpha}=\alpha\}$$
    is club, consequently, by condition \ref{4 de restrictionclubnotstationary} of Lemma \ref{restrictionclubnotstationary}, $T\restriction C\in (NS^T)^*$. For every $t\in T$ let $D_t\subseteq t\down$ be as follows:
     \[
    D_t=\left\{
    	       \begin{array}{ll}
                    A_\alpha\cap t\down      & \mathrm{if\ } t\in T\restriction C \\
                    \emptyset      & \mathrm{if\ } t\not\in T\restriction C.\\
    	       \end{array}
	     \right.
    \]

    \begin{claim}
       $(D_t\mid t\in T)$ is a $\diamT$-sequence. 
    \end{claim}
    \textit{Proof of the Claim.}
    Let $X\subseteq T$. As $T=\w_1$ then $D=\{\alpha\in\w_1\mid X\cap\alpha=A_\alpha\}$ is stationary in $\w_1$, so it is $E=D\cap C$ and thus $T\restriction E\not\in NS^T$. Let any $t\in T\restriction E$ and let $\alpha=\htop(t)$. Now we know that $X\cap T_{<\alpha}=X\cap\alpha=A_\alpha$, thus: $$X\cap t\down=(X\cap\alpha)\cap t\down=A_\alpha\cap t\down=D_t.$$
    This proves that $T\restriction E\subseteq \{t\in T\mid X\cap t\down=D_T\}$ so this last set is not in $NS^T$.\hfill\ensuremath{_{Claim}\square}    
\end{proof}


\begin{mydef}\label{def de almost Suslin}
    \cite{devlinshelahsuslintrees} Let $T$ be an $\omega_1$-tree. We say that $T$ is \textit{almost-Suslin} if, for every antichain $A\subseteq T$, the set $\widehat{A} = \{ \htop(t) \mid t \in A \}$ is nonstationary in $\omega_1$.
\end{mydef}

\begin{lemma}\label{almost-Suslin satisfies the suficient condition}
    Let $T$ be an almost-Suslin tree. Then, for every $X \subseteq \omega_1$, $T \restriction X \in NS^T$ if and only if $X$ is nonstationary.
\end{lemma}

\begin{proof}
    By condition \ref{2 de restrictionclubnotstationary} of Lemma \ref{restrictionclubnotstationary}, we know that if $X \subseteq \omega_1$ is nonstationary, then $T \restriction X \in NS^T$. Thus, it is enough to prove the converse. Suppose that $T \restriction X \in NS^T$, which means that there is a collection $\{A_t \mid t \in T\}$ of special sets of $T$ such that $T \restriction X = \nabla_{t \in T} A_t$. Hence, for every $t\in T$,  $A_t$ is a countable union of antichains, and since $T$ is almost-Suslin, every antichain is such that the set of its heights is nonstationary, meaning that $\widehat{A_t}$ is itself nonstationary.

    Now, for every $s \in T \restriction X$, let $t(s) \in T$ be a node of minimal height such that $s \in A_t \cap t\up$. For every $\alpha \in X$, let $g(\alpha) = \min\{\htop(t(s))\mid s \in T_\alpha\}$. Clearly, $g: X \to \omega_1$ is a regressive function. The proof will be complete with the following claim:
    
    \begin{claim}
        $g^{-1}(\{\beta\})$ is nonstationary for every $\beta \in \omega_1$.
    \end{claim}
    
    \textit{Proof of the Claim.} Suppose that there exist $\beta \in \omega_1$ and a stationary set $Y \subseteq X$ such that $g(\alpha) = \beta$ for all $\alpha \in Y$. This implies that for every $\alpha \in Y$, there exist $s(\alpha) \in T_\alpha$ and $t(\alpha) \in T_\beta$ such that $s(\alpha) \in A_{t(\alpha)}$. If we define $R:=\bigcup_{t \in T_\beta} A_t$, then we have
    \[
    \{s(\alpha) \mid \alpha \in Y\} \subseteq R,
    \]
    and consequently $Y \subseteq \widehat{R} = \bigcup_{t \in T_\beta} \widehat{A_t}$.
    However, $\widehat{R}$ is nonstationary (being a countable union of nonstationary sets) and contains the stationary set $Y$, which is a contradiction. \hfill \ensuremath{_{\text{Claim}} \square}
\end{proof}

\begin{theorem}\label{diamT sii diam for almost-suslin}
    If $T$ is an almost-Suslin tree, then $\diamT$ holds if and only if $\diam$ holds.
\end{theorem}

\begin{proof}
    $\diamT$ implies $\diam$ by Theorem \ref{diamT implies diam}, and $\diam$ implies $\diamT$ by Lemmas \ref{suficient condition of T for diam implying diamT} and \ref{almost-Suslin satisfies the suficient condition}.
\end{proof}

In particular:

\begin{corollary}\label{diam iff diamT if T is Suslin}
    If $T$ is a Suslin tree, then $\diamT$ holds if and only if $\diam$ holds.
\end{corollary}

The trees that arise in the study of walks on ordinals (see \cite{walksonordinals}) are important examples of \mbox{$\omega_1$-trees}. It turns out that the trees $T(\rho_0)$ and $T(\rho_2)$ are always special (see \cite{walksonordinals} and \cite{peng2013characterization}). However, the situation for $\rho_1$ is different: $T(\rho_1)$ can be nonspecial and can even be almost-Suslin. There are two different ways to see this.

Todorcevic showed that given a ladder system $\vecC$ (see Definition \ref{def de ladder system}) and $r\in([\omega]^{<\omega})^\omega$, one can associate another sequence $\vecC^r$ which, for suitable $r$, is again a ladder system. He proved that if $r$ is a Cohen real then the $\rho_1$-tree associated to $\vecC^r$ is almost-Suslin (see \cite[Theorem 15]{stevorepresentingtrees}).

A second approach uses the combinatorial principle $\bigstar_1$, which asserts the existence of a ladder system $\vecC=(C_\alpha\mid \alpha\in\mathrm{LIM}(\omega_1))$ with the property that for every stationary $S\subseteq\omega_1$ there exist ordinals $\alpha,\beta\in S\cap\mathrm{LIM}(\omega_1)$ such that $C_\beta\cap\alpha\sqsubseteq C_\alpha$ and $\alpha\in C_\beta$. Hrušák and Martínez Ranero proved that if $\vecC$ is a $\bigstar_1$-sequence then the $\rho_1$-tree associated to $\vecC$ is almost-Suslin (see \cite[Theorem 4]{hrusakranero}).

Therefore it is consistent that $T({\rho_1})$ is almost-Suslin, in particular, by Theorem \ref{diam iff diamT if T is Suslin}, we obtain:

\begin{corollary}\label{diamT is consistent for T=rho1}
    It is consistent that $\diam_{T({\rho_1})}$ holds.
\end{corollary}

Since coherent trees tend to be $\leq$-comparable, the preceding result naturally leads to the question of what can be said about $\diamondsuit_T$ for this class of trees.

\section{\texorpdfstring{$\diam^*$ is stronger than ``$\diam_T$ holds for every nonspecial Aronszajn tree''}{Diamond star is stronger that the conjunction of all diamond for nonspecial Aronszajn trees}}\label{section diamond star es mas fuerte que diamT}
We already know that for every nonspecial, well-pruned $\omega_1$-tree, $\diam^*$ implies $\diam_T$ (see Corollary \ref{diamT is consistent}). An interesting question, then, is whether it is consistent with $\neg\diam^*$ that $\diam_T$ holds for every well-pruned Aronszajn tree $T$. The first part of this section is devoted to providing a positive answer to this question. 

\begin{theorem}\label{teorema que separa diam* de diamT}
    It is consistent with $\neg\diam^*$ that there exist nonspecial, well-pruned Aronszajn trees and that $\diam_T$ holds for all such trees.
\end{theorem}

The proof of Theorem \ref{teorema que separa diam* de diamT} proceeds via a countable-support iteration of length $\omega_2$ over a model of $\mathsf{CH}$. We alternate forcing with $2^{<\omega_1}$—which destroys $\diam^*$-sequences—and the Jensen--Kurepa forcing, which introduces $\diam^+$-sequences. Although these two forcings have opposing effects, both are $\sigma$-closed; therefore, by Baumgartner's theorem (see Theorem \ref{Baumgartner}), under $\mathsf{CH}$ the iteration satisfies the $\omega_2$-chain condition. In particular, every nonspecial $\omega_1$-tree $T$ in $V[G]$ appears at some intermediate stage of the iteration. After that stage we add a $\diam^+$-sequence, which implies $\diam_T$ at that step. Moreover, the remainder of the iteration is $\sigma$-closed, so it neither specializes trees nor destroys $\diam_T$-sequences (see Lemmas \ref{sigma closed does not specialize w1-trees} and \ref{sigma closed does not kill diamT}); hence $\diam_T$ holds in the final model $V[G]$. On the other hand, also by the $\omega_2$-c.c., every $\diam^*$-candidate $\vecD$ in $V[G]$ already appears at some intermediate stage; after that stage we force with $2^{<\omega_1}$, ensuring that $\vecD$ is no longer a $\diam^*$-sequence. Since the remainder of the iteration is proper, $\vecD$ remains non-$\diam^*$ in $V[G]$.

Thus we proceed by iterating $\sigma$-closed forcing notions over a model of $\mathsf{CH}$. By Lemma \ref{sigma closed does not specialize w1-trees} such forcings do not specialize $\omega_1$-trees; we now show that they also preserve $\diam_T$-sequences.

\begin{lemma}\label{sigma closed does not kill diamT}
    Let $T$ be a nonspecial $\omega_1$-tree and $(D_t)_{t \in T}$ a $\diamondsuit_T$-sequence. If $\mathbb{P}$ is a $\sigma$-closed forcing, then $\mathbb{P}\forces``(D_t)_{t \in T} \text{ is a } \diamondsuit_T\text{-sequence}"$.
\end{lemma}

\begin{proof}
    Let $p \in \mathbb{P}$ and $\dot{X}, \dot{B}$ be such that $p \Vdash "\dot{X} \subseteq \check{T} \wedge \dot{B} \in NS^T"$. We aim to show that there exists $q \leq p$ such that $q \Vdash "\exists t \in T \setminus \dot{B} (\dot{X} \cap \check{t}\down = \check{D_t})"$.

    Now, $\dot{B} = \nabla_{s \in T} \dot{A_s}$, where each $\dot{A_s}$ is a $\mathbb{P}$-name for a special subset of $T$, there exists a $\mathbb{P}$-name $\dot{f_s}$ for a specializing function of $\dot{A_s}$, i.e., $\dot{f_s}:\dot{A_s} \to \omega$.

    For every $s, t \in T$, define 
    \[
    D(s, t) = \{q \in \mathbb{P} \mid q \text{ decides if } t \in \dot{A_s}, \text{ and if } q \forces ``t \in \dot{A_s}", \text{ then } q \text{ decides } \dot{f_s}(t)\}.
    \]
    Note that $D(s, t)$ is dense in $\mathbb{P}$ for every $(s, t) \in T^2$. 

    For every $s \in T$, let $E(s)$ denote the set $\{q \in \mathbb{P} \mid q \text{ decides } s \in \dot{X}\}$, which is also dense in $\mathbb{P}$. 

    Since $\mathbb{P}$ is $\sigma$-closed, there exists a filter $G \subseteq \mathbb{P}$ such that $G \cap D(s, t) \neq \emptyset \neq G \cap E(s)$ for every $(s, t) \in T^2$.

    For each $\P$-name $\dot{x}$, we can evaluate it with respect to the filter $G$, namely $\dot{x}[G]$, although $G$ need not be generic. Thus, for every $s\in T$ consider $\dot{A}_s[G]$ and $\dot{f}_s[G]$, both of which belong to $V$ since $G\in V$. Now set
$$
\hat{B} := \nabla_{s\in T}\bigl(\dot{A}_s[G]\bigr),
$$
and note that $\hat{B}\in NS^T$ because each $\dot{A}_s[G]$ is special, as witnessed by $\dot{f}_s[G]$.

    Next, construct a decreasing sequence $(q_\alpha)_{\alpha \in \omega_1}$ such that $q_0 \leq p$ and, for every $\alpha \in \omega_1$, we have 
    \[
    q_\alpha \Vdash``\left(\dot{B} \cap T_{<\alpha} = \hat{B} \cap T_{<\alpha}\right) \wedge \left(\dot{X} \cap T_{<\alpha} = \dot{X}[G] \cap T_{<\alpha}\right)".
    \]

    Since $\hat{B} \in NS^T$, there exists $t \in T$ such that $t \not\in \hat{B}$ and $\dot{X}[G] \cap t\down = D_t$. Let $\alpha = \htop(t) + 1$, and note that $q_\alpha$ satisfies:
    \begin{itemize}
        \item $q_\alpha \forces ``t \not\in \hat{B}"$, and thus $q_\alpha \forces ``t \not\in \dot{B}"$ because $q_\alpha \forces ``\dot{B} \cap T_{<\alpha} = \hat{B} \cap T_{<\alpha}"$.
        \item $q_\alpha \forces ``\dot{X} \cap t\down = D_t"$, since $
        q_\alpha \forces ``\dot{X} \cap T_{<\alpha} = \dot{X}[G] \cap T_{<\alpha}"$
        and $q_\alpha \forces ``\dot{X}[G] \cap t\down = D_t"$ (this last is because $\dot{X}[G] \cap t\down = D_t$ holds in $V$).
    \end{itemize}

    Hence, $q_\alpha \forces ``t \not\in \dot{B} \wedge (\dot{X} \cap t\down = D_t)"$, which completes the proof.
\end{proof}

If $\vecD=(\cD_\alpha)_{\alpha\in\w_1}$ is such that $\cD_\alpha\subseteq [\alpha]^{\leq \w}$ for every $\alpha\in\w_1$, then let us say that $\vecD$ is a \mbox{\textit{$\diam^*$-candidate}}.

The proof of the following result, which appears as Exercise (J5) in \cite[p.~300]{kunenbook}, follows essentially the same argument given in \cite[p.~387]{silverkurepaconjecture}, and it also appears in Devlin's proof that $\diam$ does not imply $\diam^*$ \cite[Theorem~3.2]{devlinvariations}.

\begin{proposition}\label{2w1 mata candidatos}
Let $\vecD$ be a $\diam^*$-candidate and $\P = 2^{<\omega_1}$. Then 
$\P \forces \text{``}\vecD\text{ is not a }\diam^*\text{-sequence.''}$
\end{proposition}

The following two forcing notions will play an important role in the proof of Theorem \ref{teorema que separa diam* de diamT}.

The \textit{Jech-Suslin forcing} $\mathbb{JS}$ is defined as: $p\in\mathbb{JS}$ if and only if $p=\emptyset$ or $p$ is a subtree of $2^{\w_1}$ such that:
\begin{enumerate}[label=(\roman*)]
    \item\label{1 de jech forcing} $\htt(p)=\alpha+1$ for some $\alpha\in\mathrm{LIM}(\w_1)$,
    \item\label{2 de jech forcing} For all $s\in p$ such that $\dom(s)<\alpha$ we have that $s\conc 0,s\conc 1\in p$,
    \item\label{3 de jech forcing} For all $\xi\leq\alpha$ we have $|p_\xi|\leq\w$ and
    \item\label{4 de jech forcing} For all $s\in p$ there exists $t\in p_\alpha$ such that $s\subseteq t$.
\end{enumerate}
And $p\leq q$ if $q=\{s\in p\mid \htop(s)<\htt(q)\}$, that is, if $q=p\restriction \htt(q)$. The reason we call $\mathbb{JS}$ the \textit{Jensen-Suslin forcing} is that Jech introduced it to force the existence of a Suslin-tree (see \cite[Exercices H11, p. 248]{kunenbook})

Note that if $p,q$ are two compatible elements in $\mathbb{JS}$ then $p\leq q$ or $q\leq p$. In fact, if $\htt(q)\leq\htt(p)$ and $r\leq p,q$ then:
$$q=r\restriction \htt(q)=(r\restriction\htt(p))\restriction \htt(q)=p\restriction\htt(q)$$
and thus $p\leq q$. It is also easy to see that $|\mathbb{JS}|=2^\w$, so
if $\mathsf{CH}$ holds, then $|\mathbb{JS}|=\w_1$.


The \textit{Jensen-Kurepa forcing} $\mathbb{JK}$ is defined as follows:
$(p,I)\in \mathbb{JK}$ if the following hold:
\begin{enumerate}
    \item $p\in\mathbb{JS}$ and $I\subseteq 2^{\w_1}$ such that $|I|\leq \w$,
    \item $f\restriction(\htt(p)-1)\in p$ for all $f\in I$.
\end{enumerate}
and $(p,I)\leq (q,J)$ if $p\leq q$ in $\mathbb{JS}$ and $I\supseteq J$. Analogously to the case of $\mathbb{JS}$, the reason we call $\mathbb{JK}$ \textit{Jensen-Kurepa} is that Jensen introduced it to force the existence of a Kurepa-tree under $\mathsf{CH}$ (see \cite[Exercices H18 and H19, p. 249]{kunenbook})

Recall that a subset $R$ of a forcing $\P$ is \emph{centered} if every finite subset of $R$ has a lower bound in $\P$, and it is \emph{linked} if every pair of elements of $R$ has a lower bound in $\P$ (in both cases the bound need not lie in $R$). We say that $\P$ is \emph{$\omega_1$-centered} (respectively, \emph{$\omega_1$-linked}) if it can be expressed as the union of $\omega_1$ many centered (respectively, linked) subsets.

Note that if we have two conditions in $\mathbb{JK}$ with the same stem, namely $(p,I)$ and $(p,J)$, then $(p,I\cup J)$ is a common lower bound. This fact, together with the observation that $|\mathbb{JS}| = 2^\omega$, shows that if $\mathsf{CH}$ holds, then $\mathbb{JK}$ is $\omega_1$-centered; in particular, it is $\omega_1$-linked. Another important consequence of $\mathsf{CH}$ about $\mathbb{JK}$ is the following:

\begin{proposition}\label{diam+ se cumple cuando fuerzas con JC}\cite[Exercice H20, p. 249 ]{kunenbook}
    If $\mathsf{CH}$ holds, then $\mathbb{JK}\forces``\diam^+"$. 
\end{proposition}

On the other hand, an important property of $\mathbb{JK}$ is that it is $\sigma$-closed. In fact, let $(p_n, I_n)_{n\in\omega}$ be a decreasing sequence in $\mathbb{JK}$, and for every $n\in\omega$ let $\alpha_n\in\mathrm{LIM}(\omega_1)$ be the unique ordinal such that $\alpha_n+1=\mathrm{ht}(p_n)$; note that $\alpha_n\leq\alpha_{n+1}$ for all $n\in\omega$.

If the sequence $(\alpha_n)_{n\in\omega}$ is eventually constant, say with constant value $\alpha_m$, then $p_n=p_m$ for every $n\geq m$, and consequently $(p_m, \bigcup_{n\in\omega} I_n)$ is a lower bound for $(p_n, I_n)_{n\in\omega}$. Now, suppose that $(\alpha_n)_{n\in\omega}$ is not eventually constant and let $\alpha=\lim_{n\in\omega} \alpha_n$. Consider the set
$$
S=\{ s\in 2^{<\omega_1} \mid (\exists m\in\omega)(s\in p_m) \wedge (\forall f\in\bigcup_{n\in\omega} I_n)(s\not\subseteq f) \}.
$$
It is clear that $|S|\leq\omega$, and by condition \ref{4 de jech forcing} of the definition of Jech-Suslin forcing, for every $s\in S$ there exists a (cofinal) branch $g\in\left[\bigcup_{n\in\omega} p_n\right]$ such that $s\subseteq g$. For each $s\in S$, fix one such branch $g_s$. Now, define
$$
r=\Bigl(\{ f\restriction\beta \mid (\exists n\in\omega)(f\in I_n) \wedge (\beta\leq\alpha) \}\Bigr) \cup \Bigl(\bigcup_{n\in\omega} p_n\Bigr) \cup \{ g_s \mid s\in S \}.
$$
We claim that $(r,\bigcup_{n\in\omega} I_n)$ is in $\mathbb{JK}$.

It clearly suffices to show that $r\in\mathbb{JS}$. Conditions \ref{1 de jech forcing} and \ref{4 de jech forcing} hold by construction, while condition \ref{3 de jech forcing} also holds because $|r_\alpha|\leq\omega$ and $|r_\xi|\leq|r_\alpha|$ by condition \ref{4 de jech forcing}. To verify condition \ref{2 de jech forcing}, note that for every $s\in r$ with $\htop(s)<\alpha$, there exists some $n\in\omega$ such that $s\in p_n$. Indeed, the only nontrivial case occurs when $s=f\restriction\beta$ for some $\beta<\alpha$ and some $f\in I_n$ for some $n\in\omega$. Since $(\alpha_n)_{n\in\omega}$ is (eventually) increasing, there exists $m>n$ such that $\alpha_m>\beta$, and because $I_n\subseteq I_m$, we have $f\restriction\alpha_m\in p_m$. Consequently, $s=f\restriction\beta\in p_m$, and hence both $s\conc 0$ and $s\conc 1$ belong to $p_m\subseteq r$, which completes the verification. The last argument also shows that $r\restriction\alpha_{m}+1=p_m$, so we are done.
 
An important definition for proving Theorem \ref{teorema que separa diam* de diamT} is the following:
\begin{mydef}
    A forcing $\P$ is called \textit{well-met} if every two compatible conditions have an infimum in $\P$.
\end{mydef}

It is clear that $2^{<\omega_1}$ is well-met since if $p, q \in 2^{<\omega_1}$ then either $p \leq q$ or $q \leq p$. It is also true that $\mathbb{JK}$ is well-met. To see this, suppose that $(p,I)$ and $(q,J)$ are compatible. In particular, as mentioned before, we must have either $p \leq q$ or $q \leq p$. Without loss of generality, assume $p \leq q$. We claim that $(p, I \cup J)$ is a condition in $\mathbb{JK}$. 

Indeed, since $(p,I)$ and $(q,J)$ are compatible in $\mathbb{JK}$, there exists some $(r,K) \leq (p,I), (q,J)$. In particular, $J \subseteq K$, and therefore for every $f \in J$ we have 
$$
f \restriction (\operatorname{ht}(r)-1) \in r.
$$ 
But since $p = r \restriction \operatorname{ht}(p)$, it follows that 
$$
f \restriction (\operatorname{ht}(p)-1) \in p.
$$ 
Thus, $(p, I \cup J)$ is a condition. Moreover, $(p, I \cup J)$ is the infimum of $(p,I)$ and $(q,J)$, since if $(r',K)$ is a lower bound for both, then necessarily $r' \leq p$ and $K \supseteq I \cup J$.

The following theorem by Baumgartner \cite{baumgartneriteraded} (which also appears as \cite[Lemma V.5.14, p. 360]{kunennuevo}) is the final piece we need to prove Theorem \ref{teorema que separa diam* de diamT}:

\begin{theorem}\label{Baumgartner}
    Let $\P_\alpha = \langle \P_\beta, \dot{\bQ}_\beta \mid \beta \in \alpha \rangle$ be a countable support iteration of forcing notions that are $\sigma$-closed, well-met and $\omega_1$-linked. If $\mathsf{CH}$ holds, then $\P_\alpha$ satisfies the $\omega_2$-chain condition.
\end{theorem}

Theorem~\ref{Baumgartner} is central to the study of the forcing axiom $\mathsf{BACH}$, which yields significant topological and set-theoretic consequences (see \cite{tall}).

\textit{Proof of Theorem \ref{teorema que separa diam* de diamT}.} Start with a model of $\mathsf{GCH}$ and consider $\P = \langle \P_\beta, \dot{\bQ}_\beta \mid \beta \in \omega_2 \rangle$ a countable support iteration of forcing notions such that at even stages we use $2^{<\omega_1}$ and at odd stages we use $\mathbb{JK}$. That is, if $\beta \in \omega_2$ is even, then
\[
\P_\beta \Vdash ``\dot{\bQ}_\beta = 2^{<\omega_1}"
\]
and if $\beta$ is odd, then
\[
\P_\beta \Vdash ``\dot{\bQ}_\beta = \mathbb{JK}".
\]

Now, let $G$ be $\P$-generic filter. By Theorem \ref{Baumgartner}, assuming $\mathsf{CH}$, the iteration $\P$ has the $\omega_2$-cc.  In particular, any $\P$-name $\dot{\vecD}$ for a $\diam^*$-candidate is decided at some stage $\gamma$.  Since at that stage or the next we force with $2^{<\omega_1}$, Proposition \ref{2w1 mata candidatos} implies that $\P_{\gamma+1}\Vdash ``\dot{\vecD}\text{ is not a }\diam^*\text{-sequence}."$. 
Hence in $V[G\restriction{\gamma+1}]$ there is $X\subseteq\omega_1$ such that $\{\alpha<\omega_1 \mid X\cap\alpha\in\dot{\mathcal D}_\alpha[G\restriction{\gamma+1}]\}$
fails to contain a club.  Equivalently, $\{\alpha<\omega_1 \mid X\cap\alpha\notin\dot{\mathcal D}_\alpha[G\restriction{\gamma+1}]\}$ is stationary.  Since each iterand in the remainder of the iteration is proper, the tail forcing is proper and thus this stationary set remains stationary in $V[G]$.  Therefore $\P\Vdash ``\dot{\vecD}\text{ is not a }\diam^*\text{-sequence}"$. Finally, because the argument applies to every $\P$-name $\dot{\vec D}$, we conclude
$$
\P\Vdash ``\neg\diam^*." 
$$

On the other hand, since $\P$ has the $\omega_2$-cc, any $\P$-name $\dot T$ for an $\omega_1$-tree in $V[G]$ is decided at some stage $\eta$.  By Lemma \ref{sigma closed does not specialize w1-trees}, if $\dot T[G\restriction\eta]$ is nonspecial then it remains nonspecial in $V[G]$; in particular, every nonspecial Aronszajn tree from $V$ stays nonspecial in $V[G]$.  Moreover, because we force cofinally often with $\mathbb{JK}$ (which adds $\diam^+$ and hence $\diam^*$), Corollary \ref{diamT is consistent} shows that $\diam_T$ holds in $V[G\restriction(\eta+1)]$.  Finally, each remaining iterand is $\sigma$-closed, so the tail forcing is $\sigma$-closed; therefore, by Lemma \ref{sigma closed does not kill diamT}, $\diam_T$ continues to hold in $V[G]$.  \hfill$\qed_{\ref{teorema que separa diam* de diamT}}$


   
\section{\texorpdfstring{$\diam$ may be strictly weaker than some $\diamT$}{Diamond may be strictly weaker than diamond of some tree}}

We already know that if $T$ is a nonspecial $\w_1$-tree then $\diamT$ is strictly weaker than $\diam^*$ and is at least as strong as $\diam$; in fact, if $T$ is Suslin then $\diamT$ is equivalent to $\diam$. Thus, the remaining question is whether there exists (or may exist) a nonspecial $\w_1$-tree such that $\diam$ is strictly weaker than $\diamT$. This entire section is dedicated to give a positive answer to this.

\begin{theorem}\label{conjectura diamante mas debil que diamante T}
    It is consistent with $\diam$ that there exists a nonspecial Aronszajn tree $T$ such that $\diamT$ does not hold.
\end{theorem}

The plan to prove Theorem \ref{conjectura diamante mas debil que diamante T} is as follows:
\begin{enumerate}
    \item Start with $R\subseteq\w_1$ being bistationary (i.e., both $R$ and $\omega_1 \setminus R$ are stationary) and a model of $\diam^*{(\omega_1 \setminus R)}$ (for example $V = L$).
    \item Thanks to Theorem \ref{important theorem of Shelah}, we have a nonspecial Aronszajn tree $T$ such that $T \restriction R \in NS^T$ and $T \restriction (\omega_1 \setminus R) \notin NS^T$.
    \item Iterate the forcing from Definition \ref{club guessing forcing} to kill $\clubsuit(\omega_1\setminus R)$\footnote{See Definition \ref{def of club guessing}; for now it suffices to know that $\diam(\omega_1\setminus R)$ implies $\clubsuit(\omega_1\setminus R)$.}. This is carried out in Proposition \ref{P(S,C) mata a S como sucesion club} and Theorem \ref{diamante del complemento de S y no club en S}\footnote{These results are applied with $S = \omega_1\setminus R$.}
    
    \item Prove that this iteration of forcing preserves the following:
    \begin{enumerate}
        \item $\diam(R)$ (Theorem \ref{diamante del complemento de S y no club en S} applied to $S:=\w_1\setminus R$)),
        \item $T \restriction (\omega_1 \setminus R) \notin NS^T$ (By Theorem \ref{statclosedinmodels forcings no mata estacionarios de arboles}, Theorem \ref{iteracion de scm es scm} and Proposition \ref{P is statclosedinmodels}).
    \end{enumerate}
    \item Apply Theorem \ref{improving diamT implies diam} to argue that $\diamT$ does not hold in $V[G]$, since in this case $\diam(\w_1\setminus R)$ would hold in $V[G]$, but not even $\clubsuit{(\w_1\setminus R)}$ holds in $V[G]$.
    \item Then such a model is a model of $\diam(R)+\neg 
    \diamT$, in particular a model of $\diam+\neg \diamT$.
\end{enumerate}

We will begin by recalling an important theorem of Shelah. The theorem essentially states that, under certain hypothesis, given a bistationary set $S \subseteq \omega_1$, there exists a nonspecial Aronszajn tree $T$ such that $T \restriction S \in NS^T$. Although this theorem is originally stated in different terms, it is crucial for the proof of Theorem \ref{conjectura diamante mas debil que diamante T}, since, as previously suggested, such a tree will serve—after a forcing extension—as the witness for the theorem. As Shelah’s theorem is formulated using different terminology, we will first introduce a preliminary definition and lemma in order to properly state it.

\begin{mydef}\cite[Definition 3.3, p. 444]{properandimproper}
 Let $S \subseteq \omega_1$ be stationary, and let $T$ be an $\omega_1$-tree. We say that $T$ is \textit{$S$-st-special} if there exists a function $f$ such that:
\begin{enumerate}
    \item $\text{dom}(f) = T \restriction S \setminus \{0\}$,
    \item if $t \in T_\alpha$ for some $\alpha \in S \setminus \{0\}$, then $f(t) \in \alpha \times \omega$, and
    \item if $s, t \in \text{dom}(f)$ are such that $s < t$, then $f(s) \neq f(t)$.
\end{enumerate}
\end{mydef}

It turns out that if $T$ is $S$-st-special for some $S\subseteq\w_1$ stationary, then $T$ is Aronszajn but not Suslin \cite[Claim 3.4(1), p.~445]{properandimproper}. On the other hand, remarkably, Shelah’s previous notion of $S$-st-special can be expressed in Brodsky’s terminology:

\begin{proposition}
Let $T$ be an $\omega_1$-tree and let $S\subseteq\omega_1$ be stationary.  Then $T$ is $S$-st-special if and only if $T\restriction S\in NS^T$.
\end{proposition}

\begin{proof}
    Suppose that $T$ is $S$-st-special, and let $f = (f_0, f_1)$ be a witness of this\footnote{Here, we are using the notation $f(t) = (f_0(t), f_1(t))$ for every $t\in T\restriction S\setminus\{0\}$.}. Now, let $g\colon T \restriction (S \setminus \{0\}) \to T$ be given by $g(s) = s \restriction f_0(s)$. Clearly, $g$ is regressive. Let us now show that $g^{-1}[\{t\}]$ is special for every $t \in T$. Note that:
\[
g^{-1}[\{t\}] = \{s \in \text{dom}(f) \mid (s > t) \wedge (f_0(s) = \htop(t))\} = \bigcup_{n \in \omega} \{s \in \text{dom}(f) \mid (s > t) \wedge (f(s) = (\htop(t), n))\}.
\]
If $s_0, s_1 > t$ are such that $f(s_0) = (\htop(t), n) = f(s_1)$, then necessarily $s_0$ and $s_1$ are incomparable, which proves that $g^{-1}[\{t\}]$ is special. This proves that $T \restriction (S \setminus \{0\}) \in NS^T$ and consequently $T \restriction S \in NS^T$.

Conversely, suppose that $T \restriction S \in NS^T$. Then, there exists a regressive function $g: T \restriction S \to T$ such that $g^{-1}[\{t\}]$ is special for every $t \in T$. That is, we can write
\[
g^{-1}[\{t\}] = \bigcup_{n \in \omega} A^t_n,
\]
where each $A^t_n$ is an antichain. Now, let us define a function $f$ as follows:
\begin{enumerate}
    \item $\text{dom}(f) = T \restriction S \setminus \{0\}$,
    \item if $s \in g^{-1}[\{t\}]$, then $f(s) = (\htop(t), n)$, where $n \in \omega$ is such that $s \in A^t_n$.
\end{enumerate}
It should be clear that $f$ witnesses that $T$ is $S$-st-special.
\end{proof}

\begin{theorem}\label{important theorem of Shelah}
    \cite[Lemma 3.9, p. 448]{properandimproper} Let $S\subseteq\w_1$ be stationary and assume that $\diam^*{(\w_1\setminus S)}$ holds. Then there is a well-pruned Aronszajn tree $T$ such that $T\restriction S\in NS^T$ and for all $R\subseteq\w_1$ we have that $T\restriction R\in NS^T$ if and only if $R\setminus S$ is nonstationary. Moreover, there is no antichain $A\subseteq T$ such that $\widehat{A}\setminus S$ is stationary.  
\end{theorem}


The fundamental forcing notion that we will work with in the remainder of this article for the proof of Theorem \ref{conjectura diamante mas debil que diamante T} is the forcing designed to destroy $\clubsuit(S)$ sequences.

\begin{mydef}\label{def de ladder system}
    Let $S \subseteq \omega_1$ be a stationary set. 
    A sequence $\vecC=(C_\alpha \mid \alpha \in S \cap \mathrm{LIM}(\omega_1))$ is called a \textit{ladder system in $S$} if $C_\alpha \subseteq \alpha$ is cofinal of order type $\omega$ for every $\alpha \in S \cap \mathrm{LIM}(\omega_1)$.
\end{mydef}

\begin{mydef}\label{def of club guessing}
    Let $S \subseteq \omega_1$ be a stationary set and $\vecC=(C_\alpha \mid \alpha \in S \cap \mathrm{LIM}(\omega_1))$ a ladder system on $S$. Then $\vecC$ is called a $\clubsuit(S)$-sequence if for every $X \in [\omega_1]^{\omega_1}$, the set $\{\alpha \in S \cap \mathrm{LIM}(\omega_1) \mid C_\alpha \subseteq X \cap \alpha\}$ is stationary.
\end{mydef}

\begin{mydef}\label{club guessing forcing}
    Let $S \subseteq \omega_1$ be stationary and $\vecC=(C_\alpha \mid \alpha \in S \cap \mathrm{LIM}(\omega_1))$ be a ladder system on $S$. The forcing $\mathbb{P}(S,\vecC)$ is defined as follows: $p \in \mathbb{P}(S,\vecC)$ if and only if there is some $\alpha \in \omega_1$ such that $p: \alpha \to 2$, and for every $\beta \in (\alpha+1) \cap S \cap \mathrm{LIM}(\omega_1)$, we have $C_\beta \not\subseteq p^{-1}(\{1\})$. The order is given by $p \leq q$ if $p \supseteq q$.
\end{mydef}

The forcing $\P(S,\vec C)$ has two desirable features: it is proper and $\omega$-distributive, hence totally proper\footnote{Recall that a forcing is \emph{totally proper} if it is proper and does not add new reals.}. Indeed, by Theorem 7.13 of \cite{justintodddavid}, the poset is even completely proper and ${<}{\omega_1}$-proper via the fusion schemes (CP$_{\aleph_1}$) and (A). It also has two additional properties crucial for our goals: it is $\sigma(\omega_1\setminus S)$-closed (Definition \ref{def de simga S closed}), and it is strategically closed in models (Definition \ref{def de stat closed in models y wscm}). In fact, we will give an alternative proof—using strategic closure in models—of its total properness (see Proposition \ref{P is statclosedinmodels} and Theorem \ref{wscm implica proper y w-dist}). However, the most important feature—and the very reason we employ this forcing—is that it is designed to destroy the sequence $\vecC$ as an instance of $\clubsuit(S)$:

\begin{proposition}\label{P(S,C) mata a S como sucesion club}
     Let $S \subseteq \omega_1$ be stationary and $\vecC=(C_\alpha \mid \alpha \in S \cap \mathrm{LIM}(\omega_1))$ be a ladder system on $S$. Then $\P(S,\vecC)\Vdash``\vecC \text{ is not a } \clubsuit(S)\text{-sequence}"$.
\end{proposition}
\begin{proof}
    For every $\alpha\in\w_1$, let $D_\alpha=\{p\in\P(S,\vecC)\mid \exists\beta\exists\gamma(\alpha\leq \gamma<\beta(p:\beta\to 2\wedge p(\gamma)=1))\}$.
    \begin{claim}
        $D_\alpha$ is dense for every $\alpha\in\w_1$. 
    \end{claim}
\textit{Proof of the Claim.} Let $p\in\P(S,\vecC)$ and let $\xi=\dom(p)$. Now, let $\beta\in S\cap\LIM(\w_1)$ be such that $\beta>\alpha,\xi$ and enumerate $(\xi,\beta]$ as $\{\alpha_n\mid n\in\w\}$ such that $\alpha_0=\beta$. Now construct recursively $\{a_n\mid n\in\w\}$ and $\{b_n\mid n\in\w\}$ as follows:
    \begin{enumerate}
        \item $a_0,b_0\in C_{\alpha_0}\setminus\xi$ such that $\alpha<a_0<b_0$.
        \item $a_{n+1},b_{n+1}\in C_{\alpha_{n+1}}\setminus(\xi\cup \{a_i,b_i\mid i\leq n\}$ and $a_{n+1}<b_{n+1}$.
    \end{enumerate}
    Now let $q:\beta\to 2$ be given by:
     \[
    q(\delta)=\left\{
    	       \begin{array}{ll}
                    p(\delta)      & \mathrm{if\ } \delta\in\xi,\\
                    0      & \mathrm{if\ } \exists n\in\w (\delta=a_n),\\
                    1      & \mathrm{if\ } \exists n\in\w (\delta=b_n),\\
                    0      & \text{in other case.\ }\\
    	       \end{array}
	     \right.
    \]
    Clearly $q\in D_\alpha$ and $q\leq p$.\hfill\ensuremath{_{Claim}\square}

    Let us go to $V[G]$. Let $f_G:=\cup G:\w_1\to 2$. Then, since every $D_\alpha$ is dense, $X_G:=f_G^{-1}(\{1\})$ is cofinal in $\w_1$ (and $\omega_1^V = \omega_1^{V[G]}$, which follows from either the properness of $\P(S,\vec C)$ or its $\omega$-distributivity).
    \begin{claim}
        $C_\alpha\not\subseteq X_G$ for every $\alpha\in\w_1$.
    \end{claim}
    \textit{Proof of the Claim.}
    If $C_\alpha\subseteq X_G$ for some $\alpha\in\w_1$, then $C_\alpha\subseteq X_G\cap\alpha$ as $C_\alpha\subseteq\alpha$. Now, let $p\in G$ be such that $\dom(p)>\alpha$. Then $p$ is such that $X_G\cap\alpha\subseteq p^{-1}(\{1\})$. Thus, we have that:
\[     C_\alpha\subseteq X_G\cap\alpha\subseteq p^{-1}(\{1\}),     \]
which is impossible as $p\in\P(S,\vecC)$.
    \hfill\ensuremath{_{Claim}\square}
\end{proof}


The next step in the proof of Theorem \ref{conjectura diamante mas debil que diamante T} is to prove that $\diam(\omega_1\setminus S)$ is preserved under iterations of forcings of the form $\P(S,\vecC)$, where $\vecC$ ranges over all possible ladder systems on $S$. To this end, we begin by proving that $\diam$ implies the existence of a seemingly stronger sequence. The following is based on a result by Hrušák and the first author \cite[Lemma 4.9]{statbounding} (see also \cite{lifeinsacksmodelMICHAEL}). 

\begin{lemma}\label{analogo de 4.9}
    Let $S\subseteq\w_1$ be stationary, $\mathbb{P}$ be a forcing notion, $\kappa$ be a large enough regular cardinal and assume $V \models \diamondsuit(S)$. Then there is a sequence $\langle (M_\alpha, p_\alpha, \dot{X}_\alpha)\mid{\alpha \in S\cap\mathrm{LIM}(\omega_1)}\rangle$ such that for every $\alpha \in S\cap\mathrm{LIM}(\omega_1)$, the following holds:
\begin{enumerate}
    \item\label{a de www} $M_\alpha$ is a countable elementary submodel of $\mathsf{H}(\kappa)$ such that $\mathbb{P}, p_\alpha, \dot{X}_\alpha \in M_\alpha$.
    \item\label{b de www} $p_\alpha \in \mathbb{P}$ and $p_\alpha \Vdash``\dot{X}_\alpha\subseteq\w_1"$.
\end{enumerate}

The sequence $\langle (M_\alpha, p_\alpha, \dot{X}_\alpha)\mid{\alpha \in S\cap\mathrm{LIM}(\omega_1)}\rangle$ has the property that for every $p \in \mathbb{P}$ and $\dot{X}$ such that $p \Vdash``\dot{X}\subseteq\omega_1"$, there is a countable $N \preceq \mathsf{H}(\kappa)$ and $\alpha\in S\cap\mathrm{LIM}(\omega_1)$ such that the following conditions hold:
\begin{enumerate}[label=(\Roman*)]
    \item $\mathbb{P}, p, \dot{X} \in N$.
    \item $M_\alpha \cap \omega_1 = \alpha$.
    \item The structures $
    (N, \in, \mathbb{P}, \leq_\P, \Vdash_\P, p, \dot{X})$ and $(M_\alpha, \in, \mathbb{P}, \leq_\P, \Vdash_\P, p_\alpha, \dot{X}_\alpha)$
    are isomorphic.
\end{enumerate}
\end{lemma}

\begin{proof}
Using $\diamondsuit(S)$ we can find a sequence
\[
\langle \mathfrak{A}_\alpha = (\alpha, \triangle_\alpha, P_\alpha, \precsim_\alpha, \rightsquigarrow_\alpha, r_\alpha, Y_\alpha) \mid{\alpha \in S\cap\mathrm{LIM}(\omega_1)}\rangle
\]
such that for every structure $\mathfrak{A} = (\omega_1, \triangle, P, \precsim ,\rightsquigarrow, r, Y)$ --with $\triangle,\precsim,\rightsquigarrow\subseteq{\w_1}^2$, $P,Y\subseteq\w_1$ and $r\in\w_1$-- there are stationary many $\alpha\in S$ such that $\mathfrak{A}_\alpha$ is an elementary substructure of $\mathfrak{A}$. Given $\alpha$ a limit ordinal in $S$, in case there are a countable $M \preceq \mathsf{H}(\kappa)$, $p \in \mathbb{P}$, $\dot{X}$ such that
\begin{enumerate}[label=(\roman*)]
    \item\label{a de bbb} $\mathbb{P},\leq_\P, p, \dot{X} \in M$, 
    \item\label{b de bbb} $M\cap\w_1=\alpha$,
    \item\label{c de bbb} $p \Vdash`` \dot{X}\subseteq\omega_1"$ and
    \item\label{d de bbb} $(M, \in, \mathbb{P}, \leq_\P,\Vdash_\P, p, \dot{X})$ is isomorphic to $\mathfrak{A}_\alpha$,
\end{enumerate}
then we choose one of them and define $M_\alpha = M$, $p_\alpha = p$ and $\dot{X}_\alpha = \dot{X}$. If there is no $M$ that satisfies these properties, we just take any $(M_\alpha, p_\alpha, \dot{X}_\alpha)$ that satisfies the properties (\ref{a de www}) and (\ref{b de www}) of the Lemma. 

We will now prove that $\langle(M_\alpha, p_\alpha, \dot{D}_\alpha) \mid \alpha \in S\cap\mathrm{LIM}(\omega_1)\rangle$ has the desired properties.

Let $p \in \mathbb{P}$ and $\dot{X}$ be a $\P$-name such that $p \Vdash``\dot{X}\subseteq\omega_1"$. Recursively, we build $\{N_\alpha \mid \alpha < \omega_1\}$ a continuous $\in$-chain of countable elementary submodels of $\mathsf{H}(\kappa)$ such that $p, \dot{X}, \mathbb{P} \in N_0$. Let $N = \bigcup_{\alpha < \omega_1} N_\alpha$, since $N$ has size $\omega_1$, then we can define a structure
\[
\mathfrak{A} = (\omega_1, \triangle, P, \leq, \rightsquigarrow,r, Y)
\]
that is isomorphic to $(N, \in, \mathbb{P},\leq_\P, \Vdash, p, \dot{X})$. Let $F : \omega_1 \to N$ be an isomorphism.

It is easy to see that
\[
\{\alpha \in \mathrm{LIM}(\omega_1) \mid N_\alpha \cap \omega_1 = \alpha \wedge F[\alpha] = N_\alpha\}
\]
is a club. In this way, we can find $\alpha\in S\cap\mathrm{LIM}(\w_1)$ such that $F[\alpha] = N_\alpha$, $N_\alpha \cap \omega_1 = \alpha$ and $\mathfrak{A}_\alpha\preceq\mathfrak{A}$. Note that $N_\alpha$, $p$, and $\dot{X}$ satisfy conditions \ref{a de bbb}, \ref{b de bbb}, \ref{c de bbb}, and \ref{d de bbb} at stage $\alpha$.  In particular, condition \ref{d de bbb} holds, showing that $\mathfrak{A}_\alpha$ is isomorphic to both 
\[
(N, \in, \mathbb{P}, \le_\P, \Vdash_\P, p, \dot{X})
\quad\text{and}\quad
(M_\alpha, \in, \mathbb{P}, \le_\P, \Vdash_\P, p_\alpha, \dot{X}_\alpha).
\]
Hence those two structures are isomorphic to each other.
\end{proof}

The following definition captures what we consider to be one of the essential properties of forcings of the form $\P(S,\vecC)$, and this property ensures that the forcings possessing it preserve certain diamond principles (see Theorem \ref{S closed forcings preserva diamante(S)} and Proposition \ref{P(S,C) is S-closed}). We must clarify that this appears to be a specific case of the theory that Shelah develops in \cite[Chapter V, Section 1]{properandimproper}; however, in our opinion, that presentation is very technical and hard to follow, which is why we decided to include this material, a presentation that we believe is quite accessible.

\begin{mydef}\label{def de simga S closed}
    Let $(\P,\leq_\P)$ be a forcing, $\kappa$ a cardinal, $M\preceq\mathsf{H}(\kappa)$ with $\P\in M$ and $(p_n)_{n\in\w}$ a $\leq_\P$-decreasing sequence. We say that $(p_n)_{n\in\w}$ is $(M,\P)$-\textit{generic} if:
    \begin{enumerate}
        \item $p_n\in M$ for all $n\in\w$ and
        \item for all $D\in M$ open dense subset of $\P$, there is $n\in\w$ such that $p_n\in D$.
    \end{enumerate}
\end{mydef}

\begin{mydef}
    Let $\P$ be a forcing and $S\subseteq\w_1$ stationary. We say that $\P$ is \mbox{$\sigma(S)$-\textit{closed}} if for all $\kappa$ large enough cardinal and all $M\preceq\mathsf{H}(\kappa)$ countable such that $\P\in M$ and $M\cap\w_1\in S$, we have that every $(M,\P)$-generic sequence has a lower bound.  
\end{mydef}

\begin{theorem}\label{S closed forcings preserva diamante(S)}
    If $\P$ is $\sigma(S)$-closed and $V\models \diam(S)$, then $\forces_\P``\diam(S)"$. 
\end{theorem}
\begin{proof}
Fix a sequence $\langle (M_\alpha, p_\alpha, \dot{X}_\alpha) \mid {\alpha \in S\cap \mathrm{LIM}(\omega_1)}\rangle$ as in Lemma \ref{analogo de 4.9}. We want to define a sequence $\mathcal{D} = \{D_\alpha\subseteq\alpha  \mid \alpha \in S\}$ that will be a witness of $\diam(S)$ in $V[G]$. Let $\alpha \in S$; in case $M_\alpha \cap \omega_1 \neq \alpha$, let $D_\alpha=\emptyset$. Now, fix $\alpha\in S$ such that $M_\alpha\cap\w_1=\alpha$ and let $\{\alpha_n\mid n\in\w\}$ be an enumeration of $\alpha$.  Fix also $\{E_n\mid n\in\w\}$ an enumeration of all open dense subsets of $\P$ in $M_\alpha$. 

Now, let us construct a sequence $\{q^\alpha_n\mid n\in\w\}$ and the set $D_\alpha\subseteq\alpha$ as follows:
\begin{enumerate}
    \item $q^\alpha_0= p_\alpha$.
    \item If $q_n^\alpha$ has been defined, then:
    \begin{enumerate}
        \item If there exists $q\in\P$ such that:
        \begin{enumerate}[label=(\roman*)]
            \item\label{1 de sub x} $q\leq q^\alpha_n$,
            \item\label{3.5 de sub x} $q\in E_n$ and
            \item\label{4 de sub x} $q\forces``\alpha_n\in \dot{X_\alpha}"$.
        \end{enumerate}
    then let $q^\alpha_{n+1}$ any such $q\in M_\alpha$ and set $\alpha_n\in D_\alpha$. 
    \item If not, then let $q^\alpha_{n+1}\in M_\alpha$ be such that satisfies conditions \ref{1 de sub x} and \ref{3.5 de sub x} and set $\alpha_n\not\in D_\alpha$.    
    \end{enumerate}
\end{enumerate}
    
       Note that this process completely determines the set $D_\alpha \subseteq \alpha$, and, moreover, $D_\alpha$ is in the ground model and it satisfies the following: 
       \begin{enumerate}[label=(\Roman*)]
           \item\label{1 de alguna cosa x} $\alpha_n\in D_\alpha$ if and only if $q^\alpha_{n+1}\forces``\alpha_n\in\dot{X_\alpha}"$ and
           \item\label{2 de alguna cosa x} $\alpha_n\not\in D_\alpha$ if and only if $q^\alpha_{n+1}\forces``\alpha_n\not\in\dot{X_\alpha}"$.
       \end{enumerate}
       
       \begin{claim}
     $(D_\alpha\mid \alpha\in S)$ is a $\diam(S)$-sequence in $V[G]$.
 \end{claim}

\textit{Proof of the Claim.} Let $\dot{X}$ be a $\P$-name and $p\in\P$ such that $p\forces``\dot{X}\subseteq\w_1"$. We want to prove that there exist $q\leq p$ and $\alpha\in S\cap\mathrm{LIM}(\w_1)$ such that $q\forces``\dot{X}\cap\alpha=D_\alpha"$.

By hypothesis, there exist $\alpha\in S\cap\mathrm{LIM}(\w_1)$ and $N\preceq \mathsf{H}(\kappa)$ such that $\mathbb{P}, p, \dot{X} \in N$, $M_\alpha \cap \omega_1 = \alpha$ and $
    (N, \in, \mathbb{P},\leq_\P, \Vdash, p, \dot{X})\cong(M_\alpha, \in, \mathbb{P},\leq_\P, \Vdash, p_\alpha, \dot{X}_\alpha)$. So let $F:M_\alpha\to N$ be the (unique) isomorphism, in particular, $F(\P)=\P$, $F(p_\alpha)=p$ and $F(\dot{X_\alpha})=\dot{X}$. Note that since $F$ is an isomorphism then, for all $n\in\w$, the following conditions hold:
    \begin{enumerate}[label=(\Alph*)]
    \item\label{a de alguna cosa x} $q^\alpha_{n+1}\forces``\alpha_n\in\dot{X_\alpha}"$ if and only if $F(q^\alpha_{n+1})\forces``\alpha_n\in\dot{X}"$,
    \item\label{b de alguna cosa x} $q^\alpha_{n+1}\forces``\alpha_n\not\in\dot{X_\alpha}"$ if and only if $F(q^\alpha_{n+1})\forces``\alpha_n\not\in\dot{X}"$ and
    \item $F(q^\alpha_{n+1})\leq F(q^\alpha_n)$. 
\end{enumerate}

Now, by construction, the sequence $(q^\alpha_n)_{n\in\omega}$ is $M_\alpha$-generic, so $(F(q^\alpha_n))_{n\in\omega}$ is $N$-generic. Since $\P$ is \mbox{$\sigma(S)$-closed}, the sequence $(F(q^\alpha_n))_{n\in\omega}$ admits a lower bound $r$. Moreover, we have
$$
r \;\le\; F(q^\alpha_0) \;=\; F(p_\alpha) \;=\; p.
$$
We claim that $r\forces``D_\alpha=\dot{X}\cap\alpha"$. To see this, note that by Conditions \ref{1 de alguna cosa x} and \ref{a de alguna cosa x} we have:
\[
\alpha_n\in D_\alpha\implies q_{n+1}^\alpha\forces``\alpha_n\in\dot{X_\alpha}"\implies F(q_{n+1}^\alpha)\forces``\alpha_n\in\dot{X}"\implies r\forces``\alpha_n\in\dot{X}"
\]
and by Conditions \ref{2 de alguna cosa x} and \ref{b de alguna cosa x}:
       \[
\alpha_n\not\in D_\alpha\implies q_{n+1}^\alpha\forces``\alpha_n\not\in\dot{X_\alpha}"\implies F(q_{n+1}^\alpha)\forces``\alpha_n\not\in\dot{X}"\implies r\forces``\alpha_n\not\in\dot{X}".
\]
       \finishclaim
\end{proof}

We want to prove that the iteration theorem for $\sigma(S)$-closed forcings, that is, that a csi of $\sigma(S)$-closed forcings is $\sigma(S)$-closed. To this end, we must first establish some preliminary results:

\begin{lemma}\label{easy remark on S-closed iterations}
    Let $\P_\alpha = \langle \P_\beta,\dot{\bQ}_\beta \mid \beta\in \alpha\}$ be a csi, $\kappa$ a cardinal and $M \preceq \mathsf{H}(\kappa)$ countable. If $(p_n)_{n \in \omega}$ is a $(M,\P_\alpha)$-generic sequence and $\beta\in\alpha\cap M$, then $(p_n\restriction\beta)_{n\in\w}$ is $(M,\P_\beta)$-generic. 
\end{lemma}
\begin{proof}
    Clearly $p_{n+1}\restriction\beta \leq p_n\restriction\beta$ for all $n\in\w$ and also $p_n\restriction\beta\in M$ since $p_n,\beta\in M$. Therefore, the only thing left to prove is that for every open dense set $E\subseteq\P_\beta$ with $E\in M$, there exists some $n\in\w$ such that $p_n\restriction\beta\in E$. To see this, note that $$
\hat{E}=\{p\in\P_\alpha \mid p\restriction\beta\in E\}\subseteq\P_\alpha$$ 
is open dense and belongs to $M$, therefore $p_n\in E$ for some $n\in\w$, thus $p_n\restriction\beta\in E$. 
\end{proof}

\begin{lemma}\label{iteration lemma case 2 and 1}
    Let $\P$ be a $\sigma(S)$-closed forcing such that $\P\forces``\dot{\bQ}\text{ is } \sigma(S)\text{-closed}"$. Let $M\preceq \mathsf{H}(\kappa)$ be a countable submodel such that $\P\ast\dot{\bQ}\in M$ and $M\cap\w_1\in S$. If $((p_n,\dot{q}_n))_{n\in\w}$ is a $(M,\P\ast\dot{\bQ})$-generic sequence and $p\in \P$ is such that $p\leq p_n$ for each $n\in\w$, then there exists $\dot{q}\in\dot{\bQ}$ such that $(p,\dot{q})\leq (p_n,q_n)$ for all $n\in\w$. 
\end{lemma}
\begin{proof}
Clearly, $p \forces ``q_{n+1} \leq q_n"$ for all $n \in \w$;  we claim that $p$ actually forces something stronger:
\begin{claim}
    $p \forces ``(\dot{q}_n)_{n\in\w} \text{ is } M[\dot{G}]\text{-generic}"$, where $\dot{G}$ is a $\P$-name for the generic filter for $\P$.
\end{claim}
\textit{Proof of the Claim.} Let $\dot{D} \in M$ be a $\P$-name for an open dense subset of $\dot{\bQ}$. Consider
\[
E \;:=\; \{(r,\dot{s}) \in \P \ast \dot{\bQ} \;\mid\; r \forces ``\dot{s} \in \dot{D}"\}.
\]
It is straightforward to check that $E$ is open dense, and $E \in M$ because $\dot{D} \in M$. Hence, there is some $n \in \w$ such that $(p_n,\dot{q}_n) \in E$. Since $p \leq p_n$, we conclude $p \forces ``\dot{q}_n \in \dot{D}"$.
\finishclaim

Therefore,
\[
p \forces ``\bigl((\dot{q}_n)_{n\in\w}\text{ is } M[\dot{G}]\text{-generic}\bigr) 
\;\wedge\; \bigl(\dot{\bQ}\text{ is }\sigma(S)\text{-closed}\bigr)",
\]
and consequently,
\[
p \forces ``\exists \dot{q} \in \dot{\bQ}\,\bigl(\forall n \in \w,\;\dot{q} \leq \dot{q}_n\bigr)".
\]
Thus, $(p,q)$ is the desired condition.
\end{proof}

The following result is the analogue of the Proper Iteration Lemma \cite[Lemma 3.3H, p. 115]{properandimproper} for $\sigma(S)$-closed forcings. 

\begin{proposition}\label{iteration lemma}
Let $\P_\alpha = \langle \P_\beta,\dot{\bQ}_\beta \mid \beta\in \alpha\}$ be a csi of $\sigma(S)$-closed forcings, $\kappa$ be a large enough cardinal, $M \preceq \mathsf{H}(\kappa)$ be a countable submodel such that $M \cap \omega_1 \in S$ and let $\gamma \in M \cap \alpha$. If $(p_n)_{n \in \omega}$ is a
$(M,\P_\alpha)$-generic sequence and $p\in \P_\gamma$ is such that $p\leq p_n\restriction\gamma$ for each $n \in \omega$, then there exists $\bar{p} \in P_\alpha$ such that $\bar{p}\restriction\gamma=p$ and $\bar{p} \leq p_n$ and for all $n \in \omega$.    
\end{proposition}

\begin{proof}
Suppose that the result is true for all $\beta<\alpha$ and let $(p_n)_{n\in\w}$ be a $(M,\P_{\alpha})$ generic sequence, $\gamma\in M\cap\alpha$ and $p\in \P_\gamma$ such that $p\leq p_n\restriction\gamma$ for all $n\in\w$.

\textbf{Case $\alpha=\beta+1$:} If $\gamma=\beta$, then the result follows directly from Lemma \ref{iteration lemma case 2 and 1}, so assume that $\gamma<\beta$. 
By Lemma \ref{easy remark on S-closed iterations} the sequence $(p_n\restriction\beta)_{n\in\w}$ is $(M,\P_{\beta})$-generic and thus, by the inductive hypothesis, there is $\bar{p}\in\P_{\beta}$ such that:
\begin{enumerate}
    \item $\bar{p}\restriction\gamma=p$ and
    \item $\bar{p}\leq p_n\restriction\beta$ for all $n\in\w$. 
\end{enumerate}
Now, since $\P_\alpha=\P_{\beta}\ast\dot{\bQ}_\alpha$ and each $p_n$ is of the form $(p_n\restriction\beta,\dot{q_n})$ where $\dot{q_n}\in\dot{\bQ}_\alpha$, we can apply Lemma \ref{iteration lemma case 2 and 1} to this sequence and condition $\bar{p}$ and obtain a condition $\hat{p}\in\P_\alpha$ such that:
\begin{enumerate}[label=(\roman*)]
    \item\label{1 de aaa} $\hat{p}\restriction\beta=\bar{p}$ and
    \item\label{2 de aaa} $\hat{p}\leq p_n$ for all $n\in\w$.
\end{enumerate}
In particular, condition \ref{1 de aaa}, implies that $\hat{p}\restriction\gamma=\bar{p}\restriction\gamma=p$.

\textbf{Case $\alpha$ is a limit ordinal:} Call $\beta:=\bigcup(\alpha\cap M)$ and let $(\alpha_n)_{n\in\w}\subseteq M$ increasing and cofinal in $\beta$ such that $\alpha_0=\gamma$. 
By recursion we will construct a sequence $(q_m)_{m\in\w}$ such that $q_0=p$ and for all $m\in\w$:
    \begin{enumerate}
        \item $q_m\in \P_{\alpha_m}$,
        \item $q_m\leq p_n\restriction\alpha_m$ for all $n\in\w$ and
        \item $q_{m+1}\restriction\alpha_m=q_m$. 
    \end{enumerate}
    Clearly $q_0$ satisfies the three conditions, so suppose that we have already constructed $q_m\in \P_{\alpha_m}$. By Lemma \ref{easy remark on S-closed iterations}, we know that $(p_n\restriction \alpha_{m+1})_{n\in\w}$ is $(M,\P_{\alpha_{m+1}})$ -generic, thus, by the inductive hypothesis there exists $q_{m+1}\in \P_{\alpha_{m+1}}$ such that the three conditions hold. Finally, let $q=\bigcup_{m\in\w}q_m$, which clearly is in $\P_\alpha$, it is a lower bound for $(p_n)_{n\in\w}$ and $q\restriction\gamma=q_0=p$. 
\end{proof} 

Applying Proposition \ref{iteration lemma} to $\gamma=0$ and $p$ the trivial condition in $\P_0=\{1\}$ we get the desired result:

\begin{theorem}\label{csi pf S-closed is S-closed}
    Let $S\subseteq\w_1$ be stationary and $\P_\alpha=\langle \P_\beta, \dot{\bQ}_\beta \mid \beta\in\alpha \rangle$ be a countable support iteration of $\sigma(S)$-closed forcings, then $\P_\alpha$ is $\sigma(S)$-closed. 
\end{theorem}

As mentioned earlier, being $\sigma(S)$-closed is, for our purposes, an important property of the forcing $\P(S,\vecC)$.

\begin{proposition}\label{P(S,C) is S-closed}
    $\P(S,\vecC)$ is $\sigma(\w_1\setminus S)$-closed.
\end{proposition}
\begin{proof}
    Let $\kappa$ be a sufficiently large cardinal and let $M \preceq \mathsf{H}(\kappa)$ be a countable elementary submodel with $M\cap\w_1=\alpha\in\w_1\setminus S$ and let $(p_n)_{n\in\w}$ be a $M$-generic sequence. To obtain the result, it suffices to prove the following:
    
    \begin{claim}
        $p:=\bigcup_{n\in\w}p_n:\alpha\to 2$ and $p\in\P(S,\vecC)$.
    \end{claim}
    
    \textit{Proof of the Claim.} Since each $p_n\in M$, we have $\dom(p_n)<\alpha$ for every $n\in\w$. Also, for every $n,m\in\w$, the conditions $p_n$ and $p_m$ are compatible, so $p$ is a function with 
    $$\dom(p)=\bigcup_{n\in\w}\dom(p_n)\leq\alpha.$$ 
    On the other hand, for any $\beta<\alpha$, the set $D_\beta=\{q\in\P(S,\vecC)\mid\dom(q)>\beta\}$ is open, dense, and belongs to $M$ since $\beta,\P\in M$. Hence there is some $n\in\omega$ with $p_n\in D_\beta$, so $\dom(p)\ge\dom(p_n)>\beta$, which shows $p\colon\alpha\to2$.

To see that $p\in\P(S,\vecC)$, we must show $C_\beta\not\subseteq p^{-1}(\{1\})$ for each $\beta\in S\cap\mathrm{LIM}(\omega_1)\cap(\alpha+1)$. Fix such a $\beta$. Since $\beta\in S$ and $\alpha\notin S$, we have $\beta<\alpha$, so there is $n\in\omega$ with $\beta<\dom(p_n)$. If $C_\beta\subseteq p^{-1}(\{1\})$, then $C_\beta\subseteq p_n^{-1}(\{1\})$, contradicting $p_n\in\P(S,\vecC)$. Hence $p\in\P(S,\vecC)$.\finishclaim
\end{proof}


The following theorem brings together some key properties of forcings of the form $\P(S,\vec C)$ and their iterations—in particular, their effects on $\diam(S)$ and $\diam(\omega_1\setminus S)$—and will be crucial for our subsequent applications.

\begin{theorem}\label{diamante del complemento de S y no club en S}
   Assume $V\models\mathsf{GCH}$.  Then there exists a forcing notion $\P$ such that $\P\Vdash`` \neg\clubsuit(S)"$.  Moreover, if $V\models\diam(\omega_1\setminus S)$, then $\P\Vdash \diam(\omega_1\setminus S)$.
\end{theorem}
\begin{proof}
    Let $\P = \langle \P_\beta, \dot{\bQ}_\beta \mid \beta \in \w_2 \rangle$ be a csi such that for each $\beta \in \w_2$, we have
\[
\P_\beta \forces ``\dot{\bQ}_\beta = \P(S,\dot{\vecC}) \text{ for some } \dot{\vecC} \text{ ladder system on } S."
\]
By a standard bookkeeping argument, we can ensure that $\P(S,\vecC)$ is eventually listed for every ladder system $\vecC$ in $S$ that appears in the intermediate models and thus, \mbox{$\forces_\P``\neg\clubsuit(S)"$}. 

Moreover, Proposition \ref{P(S,C) is S-closed} and Theorem \ref{csi pf S-closed is S-closed} imply that $\P$ is $\sigma(\omega_1\setminus S)$-closed.  Therefore, Theorem \ref{S closed forcings preserva diamante(S)} yields that $\diam(\omega_1)$ is preserved from the ground model to the extension.
\end{proof}

The final step in the proof of Theorem \ref{conjectura diamante mas debil que diamante T} is to prove that the forcing notion from Theorem \ref{diamante del complemento de S y no club en S} preserves stationary sets of trees. For this purpose, we introduce the concept of \textit{strategically-closed in models} forcings.

Let $\P$ be a forcing notion. We say that a countable model $M$ is \textit{suitable for $\P$} if $\P \in M$ and $M \preceq \mathsf{H}(\kappa)$, where $\kappa$ is a sufficiently large cardinal. 

\begin{mydef}    
Let $\P$ be a forcing and $M$ be a suitable model for $\P$. The \textit{distributivity game of $\P$ in $M$} is defined as follows:
\[
\begin{array}{c|c|c|c|c|c}
    \mathrm{I}  &  p_0 \in \P \cap M &                                     & p_1 \in \P \cap M \, (p_1 \leq q_0)  & & \ldots  \\
 \hline
  \mathrm{II}  &                    & q_0 \in \P \cap M \, (q_0 \leq p_0)  & & q_1 \in \P \cap M \, (q_1 \leq p_1)
                     & \ldots 
\end{array}
\]
and $\mathrm{II}$ wins if and only if the sequence $(p_n)_{n\in\w}$ has a lower bound in $\P$ (equivalently, if $(q_n)_{n\in\w}$ has a a lower bound in $\P$). We denote this game by $\DG(\P,M)$.
\end{mydef}

\begin{notation}
    The sequence $(p_n)_{n\in\w}$ will be called a run of the game and $(p_i)_{i\leq n}$ a partial run.
\end{notation}

\begin{mydef}\label{def de stat closed in models y wscm}
Let $\P$ be a forcing.
\begin{itemize}
    \item $\P$ is called \textit{strategically closed in models} (abbreviated as \textit{\sic}) if for every $M$ suitable model for $\P$, Player $\mathrm{II}$ has a winning strategy in $\DG(\P,M)$ 
    \item $\P$ is called \textit{weakly strategically closed in models} (abbreviated \textit{\wsic}) if for every $M$ suitable model for $\P$, Player $\mathrm{I}$ does not have a winning strategy in $\DG(\P,M)$.
\end{itemize}
\end{mydef}

As mentioned at the time of defining $\P(S,\vec C)$, its strategic closure in models is crucial for our arguments—in fact, recognizing this property inspired the very definition of strategic closure in models.

\begin{proposition}\label{P is statclosedinmodels}
    Let $S\subseteq\w_1$ be stationary and $\vecC$ be a ladder system on $S$. Then $\P(S,\vecC)$ is \sic. 
\end{proposition}
\begin{proof} 
     Call $\P=\P(S,\vec C)$. Let $M$ be a suitable model for $\P$ and set $\delta=M\cap\omega_1$. We will describe a strategy $\sigma$ for Player~II in the game $\DG(\P,M)$. 

Let $q_0$ be the first move of Player~I and choose an increasing sequence $(\alpha_n)_{n<\omega}$ such that $\alpha_0=\dom(q_0)$ and $\sup_n\alpha_n=\delta$. Now, pick $b_0\in\delta\setminus\alpha_0$ so that $b_0\in C_\delta$ if $\delta\in S$, and arbitrarily otherwise. Define
\[\sigma(q_0):\;b_0+1\to2\quad\text{by}\quad
q_0(\gamma)=\begin{cases}
  p_0(\gamma),&\gamma<\alpha_0,\\
  0,&\text{otherwise.}
\end{cases}\]
Note that $\sigma(q_0)\in M$ since $q_0,b_0\in M$. 

Suppose $(q_0,\dots,q_n)$ is a partial run of $\DG(\P,M)$ in which Player~II has followed $\sigma$. Note that there exists a condition $r\in\P$ such that $r \leq q_n$ and $\dom(r)>\alpha_n$, therefore, by elementarity and since $q_n \in M$, there is also $r' \in M$ with these properties, so set
\[\sigma(q_0,\dots,q_n)=r'\,.\]
This completes the description of $\sigma$. 

To see that $\sigma$ is a winning strategy, consider any full run $(q_n)_{n<\omega}$ where Player~II played according to $\sigma$. Let $q=\bigcup_{n<\omega}q_n$ and note that:
\begin{enumerate}  
    \item $q$ is a function,
  \item $\dom(q)=\bigcup_{n\in\w}\dom(q_n)$ and each $\dom(q_n)<\delta$ as $q_n\in M$, thus $\dom(q)\leq\delta$ and
    \item If $n\in\w$, then $\dom(q)\geq\dom(q_{n+1})\geq\dom(\sigma(q_0,\dots,q_n))>\alpha_n$/ 
\end{enumerate}
Thus $q:\delta\to2$. To see that $q\in\P$, we must show $C_\beta\not\subseteq q^{-1}(\{1\})$ for each $\beta\in S\cap\mathrm{LIM}(\omega_1)\cap(\delta+1)$. Fix such a $\beta$. If $\beta=\delta$ then we are done since $C_\delta \not\subseteq q^{-1}(\{1\})$ as $q(b_0)=(\sigma(q_0))(b_0) = 0$. If $\beta<\alpha$, there is $n\in\omega$ with $\beta<\dom(q_n)$, so, if $C_\beta\subseteq q^{-1}(\{1\})$, then $C_\beta\subseteq q_n^{-1}(\{1\})$, contradicting $q_n\in\P(S,\vecC)$. Hence $q\in\P$, so $\sigma$ is winning.  \end{proof}

An important feature of weak strategic closure in models is that it strengthens total properness: 

\begin{theorem}\label{wscm implica proper y w-dist}
    If $\P$ is a \wsicc forcing, then $\P$ is $\w$-distributive and proper. 
\end{theorem}
   \begin{proof}
To see that $\P$ is $\w$-distributive, let $p\in\P$ and let $\dot f$ be a $\P$-name with $p\Vdash``\dot f\colon \omega \longrightarrow \mathrm{ORD}"$; we must produce $q\le p$ and a ground-model function $g\colon\omega\to\mathrm{ORD}$ such that $q\Vdash``\dot f = \check g"$.

Fix a countable elementary submodel $M\prec H(\kappa)$ (for large enough $\kappa$) with $\P,\,p,\,\vecC,\,\dot f \in M$. We define a strategy $\sigma$ for Player I in the game $\DG(\P,M)$ as follows:

\begin{itemize}
  \item $\sigma(\emptyset)=p$.
  \item If Player II’s first $n+1$ moves are $(q_0,\dots,q_n)$, note that there exists a condition $r\le q_n$ such that $r$ decides the value of $\dot f(n)$. By elementarity, as $q_n,\dot{f},\P\in M$, there is such an $r'\in M$, so set \mbox{$ \sigma(q_0,\dots,q_n)=r'$}.
\end{itemize}

Since $\P$ is \wsic, Player I’s strategy $\sigma$ cannot be winning.  Hence there is a complete run of the form $(p_n)_{n<\omega}$ with
\begin{enumerate}[label=(\alph*)]
  \item $p_0=p$,
  \item each $p_{n+1}$ decides $\dot f(n)$, and
  \item the sequence $(p_n)_{n\in\w}$ has a common lower bound $q\in\P$.
\end{enumerate}
By (b) and (c), $q$ decides every $\dot f(n)$. Define $g\colon\omega\to\mathrm{ORD}$ in the ground model by
\[
g(n)=\alpha \quad\text{if}\quad q\Vdash``\dot f(n)=\alpha".
\]
Then clearly $q\Vdash``\dot f = \check g"$,
so $\P$ is $\omega$-distributive.

To see that $\P$ is proper, let $N$ be a suitable model for $\P$ and $p\in \P\cap N$. We must find a condition $q\le p$ that is $(N,\P)$-generic. Enumerate the dense open sets in $N$ as $\{D_n\mid n<\omega\}$. Define a strategy $\rho$ for Player I in the game $\DG(\P,N)$ by:
\begin{itemize}
  \item $\rho(\emptyset)=p$.
  \item If Player II’s first $n+1$ moves are $(q_0,\dots,q_n)$, note that there exists a condition $r\le q_n$ with $r\in D_n$.  By elementarity there is such an $r'\in N$, and we set $\rho(q_0,\dots,q_n)=r'$.
\end{itemize}
Since $\P$ is \wsic, $\rho$ cannot be winning.  Thus there is a full run $(p_n)_{n<\omega}$ satisfying:
\begin{enumerate}[label=(\alph*)]
  \item $p_0=p$,
  \item $p_{n+1}\in D_n\cap N$ for each $n$, and
  \item the sequence $(p_n)$ has a common lower bound $q\in\P$.
\end{enumerate}
Then $q\le p$ meets every open dense set in $N$, in particular $q$ is $(N,\P)$-generic.
\end{proof}

As we said before, from Proposition~\ref{P is statclosedinmodels} and Theorem~\ref{wscm implica proper y w-dist} we get that $\P(S,\vec C)$ is totally proper. Another important property of {\wsic} forcings is that they preserve stationary subsets of trees:

\begin{theorem}\label{statclosedinmodels forcings no mata estacionarios de arboles}
Suppose that $T$ is a nonspecial $\w_1$-tree, $X \subseteq T$ satisfies $X \notin NS^T$ and $\P$ is a \wsicc forcing. Then $\P \forces \text{``}X \notin NS^T\text{''}$. In particular $\P\forces``T\text{ is nonspecial }"$.
\end{theorem}

\begin{proof}
        Let $p\in\P$ and let $(\dot{A}_t,\dot{g}_t \mid t\in T)$ be a collection of $\P$-names such that $$p\forces \text{``}\dot{g}_t\colon\dot{A}_t\to\w \text{ is a specializing function''}$$ for every $t\in T$. We want to prove that there is some $q\le p$ and $w\in X$ such that $q\forces \text{``}w\not\in\nabla_{t\in T}\dot{A}_t\text{''}$.

For every $q\le p$, $t\in T$, and $n\in\w$, consider the set
\[
C(q,t,n)=\left\{ s\in T \,\middle|\, \Bigl(\forall z<s( q\forces \text{``}z\in\dot{A}_t \rightarrow \dot{g}_t(z)\neq n\text{''})\Bigr)
\wedge \Bigl( \exists r\le q\, \Bigl( r\forces \text{``}s\in\dot{A}_t \wedge \dot{g}_t(s)=n\text{''} \Bigr) \Bigr) \right\}.
\]

\begin{claim}
For every $q\le p$, $t\in T$, and $n\in\w$, the set $C(q,t,n)$ is an antichain.
\end{claim}

\textit{Proof of the Claim.} Suppose that $s_0,s_1\in C(q,t,n)$ with $s_0<s_1$. Since $s_1\in C(q,t,n)$, we have that $q\forces \text{``} z\in\dot{A}_t \rightarrow \dot{g}_t(z)\neq n\text{''}$, for all $z<s_1$, in particular:
\[q\forces \text{``}s_0\in\dot{A}_t \rightarrow \dot{g}_t(s_0)\neq n\text{''}.
\]
On the other hand, since $s_0\in C(q,t,n)$, there exists some $r\le q$ such that
\[
r\forces \text{``}s_0\in\dot{A}_t \land \dot{g}_t(s_0)=n\text{''}.
\]
Since $r\le q$, it follows that
\[
r\forces \text{``}(s_0\in\dot{A}_t \rightarrow \dot{g}_t(s_0)\neq n) \wedge (s_0\in\dot{A}_t \land \dot{g}_t(s_0)=n)\text{''},
\]
which is impossible.\finishclaim

Now, let $(M_\alpha \mid \alpha \in \w_1)$ be a continuous $\in$-chain of countable elementary submodels of some $\mathsf{H}(\kappa)$, for a sufficiently large $\kappa$, such that $\P$, $T$, $X$, $(\dot{A}_t,\dot{g}_t \mid t\in T)$ and $p$ belong to $M_0$. For every $t \in T$, let $\alpha(t)\in \w_1$ be minimal such that $t \in M_\alpha$, and consider
\[
B_t = \bigcup\{C(q,z,n) \mid q \in M_{\alpha(t)},\ z \leq t,\ \text{and } n \in \w\}.
\]
Note that $B_t$ is a countable union of antichains, so it is special. Let $
B := \nabla_{t \in T} B_t$.

Consider the set $
C := \{M_\alpha \cap \w_1 \mid \alpha \in \mathrm{LIM}(\w_1)\}$,  which is a club in $\w_1$. Since $T\restriction C \in (NS^T)^*$, $X \not\in NS^T$, and $B \in NS^T$, there is some 
\[
w \in \bigl((T\restriction C) \cap X\bigr) \setminus B.
\]
Consider an enumeration $\{(t_i, n_i) \mid i \in \w\}$ of the set $\{(t,n) \mid t < w \wedge n \in \w\}$ and let $\beta \in \mathrm{LIM}(\w_1)$ be such that $\htop(w) = M_\beta \cap \w_1 := \delta_\beta$. 

We will construct a strategy $\sigma$ for player $\mathrm{I}$ in the game $\DG(\P, M_\beta)$ such that $\sigma(\emptyset) = p_0:=p$ and for every $i \in \w$ we have:
\begin{equation}\label{otra ecuacion}
    p_{i+1} \forces \text{``}(w \not\in \dot{A}_{t_i}) \vee (\dot{g}_{t_i}(w) \neq n_i)\text{''}.
\end{equation}
Here, $p_{i+1}$ denotes the $(i+1)$-st move of player $\mathrm{I}$ when following the strategy $\sigma$.
 
For this suppose that in the $i$-th entry of the game player $\mathrm{II}$ plays $q_i\in\P\cap M_\beta$. We know that $t_i < \w$, i.e., $w \in t_i\up$, so, since $w \notin B$, we have $w \notin B_{t_i} \cap t_i\up$. In particular, $w \notin C(q_i,t_i,n_i)$, and thus one of the following holds:
\begin{enumerate}[label=(\Roman*)]
    \item Exists $z<w$ such that $q_i\not\forces \text{``}z \in \dot{A}_{t_i} \rightarrow \dot{g}_{t_i}(z) \neq n_i\text{''}$, or
    \item there is no $r \leq q_i$ such that $r \forces \text{``}(w \in \dot{A}_{t_i}) \land (\dot{g}_{t_i}(w) = n_i)\text{''}$.
\end{enumerate}

In other words, we obtain that one of the following holds:
\begin{enumerate}[label=(\roman*)]
    \item\label{case:i} There exist $z<w$ and and $r\leq q_i$ such that $r \forces \text{``}(z \in \dot{A}_{t_i}) \land( \dot{g}_{t_i}(z) = n_i)\text{''}$;
    \item\label{case:ii} $q_i\forces \text{``}(w \notin \dot{A}_{t_i}) \lor (\dot{g}_{t_i}(w) \neq n_i)\text{''}$.
\end{enumerate}

If case \ref{case:i} holds, then note the following:
\begin{itemize}
    \item $r$ already satisfies condition (\ref{otra ecuacion}). Indeed, if not, there exists some $r' \leq r$ such that $r' \forces \text{``}(w \in \dot{A}_{t_i}) \land (\dot{g}_{t_i}(w) = z)\text{''}$, but then $r' \forces \text{``}(z < w) \land (\dot{g}_{t_i}(z) = \dot{g}_{t_i}(w))\text{''}$, which is impossible because $r'\forces \text{``}\dot{g}_{t_i}$ is a specializing function\text{''}.
    \item By elementarity, as $z\in M_\beta$, we can take this $r \in M_\beta$ and declare $p_{i+1}:=r$.
\end{itemize}

If case \ref{case:ii} holds, then we are done by taking $p_{i+1} = q_i$. This completes the construction of the strategy $\sigma$. 

As $\P$ is a \wsicc forcing, $\sigma$ is not a winning strategy for Player $\mathrm{I}$, thus there is some run of the game where Player $\mathrm{I}$ played according to $\sigma$ and Player $\mathrm{II}$ won, that is, there exists a sequence $(p_i)_{i\in\w}$ where $p_0=p$, it satisfies the condition (\ref{otra ecuacion}) and has a lower bound $\bar{p}$.

\begin{claim}
    $\bar{p} \forces \text{``}w \notin \nabla_{t \in T} \dot{A}_t\text{''}$.
\end{claim}

\textit{Proof of the Claim.} Suppose that this is false. Then there exist $p' \leq\bar{p}$, $t < w$ and $n,i\in \w$ such that 
\begin{equation}\label{alguna ecuacion}
p' \forces \text{``}(w \in \dot{A}_t) \land (\dot{g}_{t}(w) = n)\text{''}
\end{equation}
and $(t,n) = (t_i, n_i)$. On the other hand, as $p' \leq p_{i+1}$ and $p_{i+1} \forces \text{``}(w \notin \dot{A}_{t}) \vee (\dot{g}_{t}(w) \neq n)\text{''}$, we have $$p' \forces \text{``}(w \notin \dot{A}_{t}) \vee (\dot{g}_{t}(w) \neq n)\text{''},$$ 
which contradicts (\ref{alguna ecuacion}). \finishclaim
    \end{proof}

We now have to prove the iteration theorem for \sicc forcings. As it is usual, we start by proving the case for a two-step iteration:

\begin{lemma}\label{iteracion de dos scm is scm}
    Suppose that $\P$ is a \sicc forcing and $\P\forces``\dot{\bQ}\text{ is \sic}"$. Then $\P\ast\dot{\bQ}$ is \sic.
\end{lemma}
\begin{proof}
Fix $M$ a suitable model for $\P\ast\dot{\bQ}$ and a $\P$-name $\dot{G}$ for a generic filter for $\P$. Note that:
    \begin{itemize}
        \item $M$ is also suitable for $\P$ and
        \item $\P\forces``M[\dot{G}]\text{ is a suitable model for }\dot{\bQ}[\dot{G}]"$.
    \end{itemize}
    Let $\sigma$ be a winning strategy for Player $\mathrm{II}$ in $\mathscr{DG}(\P,M)$. We will describe a strategy $\rho$ for $\mathrm{II}$ in the game $\DG(\P\ast\dot{\bQ},M)$ as follows: 
     
     Let $(p_0,\dot{q_0})\in M$ be the play of Player $\mathrm{I}$ in the first move.  We know that there exists a $\P$-name $\dot{\tau}$ such that:
     \[
     \sigma(p_0)\forces``\dot{\tau}\text{ is a winning strategy for }\mathrm{II}\text{ in the game }\DG(\dot{\bQ}[\dot{G}],M[\dot{G}])".
     \]
     Thus we have that:
     \[
    \sigma(p_0)\forces``\dot{\tau}(\dot{q_0})\in M[\dot{G}]=\{\dot{a}[G]\mid \dot{a}\in M\land(\dot{a}\text{ is a }\P\text{-name for an element of }\dot{\bQ}\}".
     \]
    In particular:
    \[
    \sigma(p_0)\forces``\exists \dot{a}\in M(\dot{a}\text{ is a }\P\text{-name for an element of }\dot{\bQ})\land(\dot{a}[G]=\dot{\tau}(\dot{q}_0)[\dot{G}])."
    \] 
    Thus we can find some $r_0\leq\sigma(p_0)$ and a $\P$-name $\dot{b}_0$ such that both are in $M$ and:
    \[
    r_0\forces``\dot{b}_0=\dot{\tau}(\dot{q}_0)."
    \]
    Hence we can take $\rho((p_0,\dot{q_0}))$ as $(r_0,\dot{b}_0)$, which by construction belongs to $M$. 

In general, if $\langle (p_i,\dot{q}_i)\mid i\leq n\rangle$ is a partial run of the game $\DG(\P\ast\dot{\bQ},M)$ where $\mathrm{II}$ played according to $\rho$ then, by the same arguments that before, we have:
\[
\sigma(p_0,\dots,p_n)\forces``\exists \dot{a}\in M(\dot{a}\text{ is a }\P\text{-name for an element of }\dot{\bQ})\land(\dot{a}[G]=\dot{\tau}(\dot{q}_0,\dots,\dot{q}_n)[\dot{G}])."
\]
 Thus we can find some $r_n\leq\sigma(p_0,\dots,p_n)$ and a $\P$-name $\dot{b}_n$ such that both are in $M$ and:
    \[
    r_n\forces``\dot{b}_n=\dot{\tau}(\dot{q}_0,\dots,\dot{q}_n)."
    \]
    Thus $\rho((p_0,\dot{q}_0),\dots, (p_n,\dot{q}_n)):=(r_n,\dot{b}_n)\in M$. 
   
    We claim that $\rho$ is a winning strategy. For this, let $\langle (p_n,\dot{q}_n)\mid n\in\w \rangle$ be a run of the game where $\mathrm{II}$ followed $\rho$ and note the following two facts:
    \begin{enumerate}
        \item for all  $n\in\w$ we have that
                \[
                r_n\forces``(\dot{q}_0\dots,\dot{q}_n)\text{ is a partial run of the game }\DG(\dot{\bQ},M[\dot{G}])\text{ where  }\mathrm{II}\text{ followed } \dot{\tau}",
                \]
    \item\label{2 de zzz} $\langle p_n\mid n\in\w\rangle$ can be seen as a run of the game $\DG(\P,M)$ where $\mathrm{II}$ followed the strategy $\sigma$. 
    \end{enumerate}    

By condition (\ref{2 de zzz}), there exists $p\in\P$ a lower bound of $(p_n)_{n\in\w}$ and this way:
    $$p\forces``\langle\dot{q}_n\mid n\in\w\rangle\text{ is a run of the game }\DG(\dot{\bQ}[\dot{G}],\dot{M}[\dot{G}]) \text{ where }\mathrm{II}\text{ followed the strategy } \dot{\tau}"$$
    and also we know that $p\forces``\dot{\tau}\text{ is a winning strategy for Player }\mathrm{II}"$, thus:
    \[
    p\forces``\exists\dot{q}\in\dot{\bQ}(\dot{q}\text{ is a lower bound of }(\dot{q}_n)_{n\in\w})".
    \]
This way $(p,\dot{q})$ is a lower bound of $\langle (p_n,\dot{q}_n)\mid n\in\w\rangle)$, which proves that $\rho$ is a winning strategy.
\end{proof}
    
Recall that if $\mathbb{B}$ is a Boolean algebra and $a,b\in\mathbb{B}$, then $a\wedge b$ denotes the infimum of $\{a,b\}$, that is, the largest lower bound of $a$ and $b$. Following this notation, if $\P_\alpha = \langle \P_\beta,\dot{\bQ}_\beta \mid \beta\in \alpha\}$ is a forcing iteration, $\beta_0<\beta_1\leq\alpha$, $p_0\in\P_{\beta_0}$, $p_1\in\P_{\beta_1}$ are such that $p_0\leq p_1\restriction\beta_0$, then we denote by $p_0\wedge p_1$ the condition $p_0\conc \left(p_1\restriction[\beta_0,\beta_1)\right)\in \P_{\beta_1}$, which can be seen as the infimum of $p_0$ and $p_1$.

\begin{theorem}\label{iteracion de scm es scm}
    Let $\P_\alpha = \langle \P_\beta,\dot{\bQ}_\beta\mid\beta\in \alpha\}$ be a countable support iteration of \sicc forcings. Then $\P_\alpha$ is \sic.
\end{theorem}
\begin{proof}
    By induction of $\alpha$. If $\alpha=\beta+1$ the result follows from Lemma \ref{iteracion de dos scm is scm}. So assume $\alpha$ is a limit ordinal, $M$ is a suitable model for $\P_\alpha$, $\dot{G}$ is a name for a generic filter for $\P_\alpha$ and $\{\alpha_m\mid m\in\w\}$ is an enumeration of $M\cap\alpha$.

    Note that $M$ is also a suitable model for $\P_{\alpha_m}$ for every $m\in\w$, so let us fix $\sigma_m$ a winning strategy for $\II$ in $\DG(P_{\alpha_m},M)$ (which exists by the inductive hypothesis) and also let us fix a $\P_{\alpha_m}$-name $\dot{\tau}_m$ such that $$\P_\alpha\forces``\dot{\tau}_m\text{ is a winning strategy for }\II \text{ in }\DG(\dot{Q}_{\alpha_m}[\dot{G}\restriction\alpha_m],M[\dot{G}\restriction\alpha_m])".$$
We want to show that there exists a strategy $\rho$ that is winning for $\II$ in $\DG(\P_\alpha,M)$. To prove this, we will begin by describing some properties that our desired strategy $\rho$ will satisfy; then, we will show that if $\rho$ has these properties, it is winning; and finally, we will demonstrate that such a $\rho$ can be constructed (which is the most technical part of the proof).

\textbf{Promises of the strategy $\rho$:}

We will describe the strategy $\rho$ for $\II$ in $\DG(\P_\alpha,M)$ by considering sequences 
$$\langle(s^m(n))_{n\in\w} \mid m\in\w\rangle,\quad \langle(q^m(n))_{n\in\w} \mid m\in\w\rangle,\quad \text{and} \quad \langle(\dot{b}^m(n))_{n\in\w} \mid m\in\w\rangle,$$
such that for every $n,m\in\w$ we have:
$$s^m(n)\in\P_\alpha,\quad q^m(n)\in\P_{\alpha_m},\quad\text{and}\quad \dot{b}^m(n)\in\dot{\bQ}_{\alpha_m}\cap M.$$
Moreover, these sequences satisfy that if $(p_0,\dots, p_m)$ is a partial run of $\DG(\P_\alpha,M)$ played according to $\rho$, that is,
$$
\begin{array}{c|c|c|c|c|c|c|c}
    \mathrm{I}  &  p_0  &                                     & p_1  & & \ldots & p_m &  \\
 \hline
  \mathrm{II}  &                    & \rho(p_0)  & & \rho(p_0,p_1) & &  & \rho(p_0,\dots,p_m)
\end{array},
$$
then the following properties hold:
\begin{enumerate}[label=(\Alph*)]
    \item For every $i\leq m$ and every $n\leq i$ we have:
        \begin{enumerate}[label=(\alph*)]
            \item\label{a chiquita de alguna} $\rho(p_0,\dots,p_i)\leq s^n(i-n)\leq p_i$, and
            \item\label{b chiquita de alguna} $\rho(p_0,\dots, p_i)\restriction \alpha_n\leq q^n(i-n)\leq s^n(i-n)\restriction\alpha_n$.
        \end{enumerate}
        \item\label{c de las q} For every $i\leq m$ and every $k<l\leq i$ we have:
        \begin{itemize}
            \item If $\alpha_k<\alpha_l$, then $q^k(i)\conc \dot{b}^k(i)\leq q^l((i+k)-l)\restriction(\alpha_k+1)$.
            \item If $\alpha_l<\alpha_k$, then $q^l(i-l)\conc \dot{b}^l(i-l) \leq q^k(i-k)\restriction(\alpha_l+1)$.
        \end{itemize}
    \item For every $n\leq m$ we have that
            $$
                \begin{array}{c|c|c|c|c|c|c|c}
                \mathrm{I}  &  s^n(0)\restriction\alpha_n  &                                     & s^n(1)\restriction\alpha_n  & & \ldots & s^n(m-n)\restriction\alpha_n &  \\
                \hline
                \mathrm{II}  &                    & q^n(0)  & & q^n(1) &\ldots &  & q^n(m-n)
            \end{array}
            $$
        is a partial run of the game $\DG(\P_{\alpha_n},M)$ such that for every $i\leq m-n$:
            \begin{itemize}
                \item $q^n(i)\leq\sigma_m\bigl(s^n(0)\restriction\alpha_n,\dots,s^n(i)\restriction\alpha_n\bigr)$, and
                \item $q^n(i)\forces\text{``}\dot{b}^n(i)=\dot{\tau}_m\bigl(s^n(0)(\alpha_n),\dots,s^n(i)(\alpha_n)\bigr)\text{''}$.
            \end{itemize}
\end{enumerate}

\textbf{Showing that $\rho$ is winning:} Note that if $(p_n)_{n\in\w}$ is a run played according to the strategy $\rho$ and $n,m\in\w$ with $n\ge m$, then by conditions \ref{a chiquita de alguna} and \ref{b chiquita de alguna} we have that:
\begin{equation}\label{ec importante}
 q^m(n-m)\le p_n\restriction\alpha_m \tag{*}   
\end{equation}

\begin{claim}
     $\rho$ is a winning strategy.
\end{claim}
\textit{Proof of the Claim.} Suppose that $(p_n)_{n\in\w}$ is a run of $\DG(\P_\alpha,M)$ in which $\II$ followed $\rho$. We will exhibit a lower bound $r$ for $(p_n)_{n\in\w}$ by constructing a sequence $\langle r_\beta \mid \beta\in\alpha\cap M\rangle$ such that for every $\beta\in\alpha\cap M$ we have:
\begin{enumerate}[label=(\roman*)]
    \item $r_\beta\in\P_\beta$,
    \item\label{2 de ggg} If $\beta=\alpha_m$, then $r_\beta$ is a lower bound of $(q^m(n))_{n\in\w}$,
    \item\label{3 de ggg} $r_\beta\forces\text{``}(\dot{b}^m(n))_{n\in\w}$ is a run of the game $\DG(\dot{\bQ}_{\alpha_m}[\dot{G}\restriction\alpha_m],M[\dot{G}\restriction\alpha_m])$, in which $\II$ followed $\dot{\tau}_m\text{''}$, and
    \item If $\beta<\gamma$ with $\beta,\gamma\in M\cap\alpha$, then $r_\beta=r_\gamma\restriction\beta$.
\end{enumerate}
Note that condition \ref{3 de ggg} is really a consequence of condition \ref{2 de ggg}. Suppose that $r_\beta$ has been constructed for every $\beta<\gamma$ in $M\cap\alpha$. 

\textbf{Case $\gamma=\beta+1$:} By condition \ref{2 de ggg} and the fact that 
$$r_\beta\forces\text{``}\dot{\tau}_m\text{ is a winning strategy for }\II\text{ in }\DG(\dot{\bQ}_{\alpha_m}[\dot{G}\restriction\alpha_m],M[\dot{G}\restriction\alpha_m])\text{''},$$
we have 
$$r_\beta\forces\text{``}\exists \dot{b}_m\in\dot{\bQ}\text{, a lower bound for }(\dot{b}^m(n))_{n\in\w}\text{''}.$$
Then we let $r_{\beta+1}=r_\beta\conc\dot{b}_m$. 

Suppose that $\beta+1=\alpha_k$. We must prove that $r_{\beta+1}$ is a lower bound for $(q^k(n))_{n\in\w}$.

\textbf{Subcase $k>m$:} By condition \ref{c de las q}, we know that for all $i\in\w$ we have
$$q^m((k-m)+i)\conc\dot{b}^m((k-m)+i)
\leq q^k(i)\restriction(\alpha_m+1)=q^k(i)\restriction \alpha_k=q^k(i).$$
On the other hand, since $r_{\beta+1}\restriction \beta=r_\beta\le q^m((k-m)+i)$ and $r_\beta\forces\text{``}\dot{b}_m\le \dot{b}^m((k-m)+i)\text{''}$, it follows that 
$$r_{\beta+1}\le q^m((k-m)+i)\conc\dot{b}^m((k-m)+i),$$
so we are done. 

\textbf{Subcase $k<m$:} In this case, by condition \ref{c de las q}, for all $i\in\w$ we have
$$q^m(i)\conc\dot{b}^m(i)\leq q^k((m-k)+i)\restriction(\alpha_m+1)=q^k((m-k)+i)\restriction(\alpha_k)=q^k((m-k)+i).$$
Again, by the same argument as before, we obtain 
$$r_{\beta+1}\le q^m(i)\conc\dot{b}^m(i),$$
so $r_{\beta+1}\leq q^k(n)$ for all $n\geq m-k$, which as $q^k(m-k)\leq q^k(j)$ for all $j\leq m-k$, proves that $r_{\beta+1}$ is a lower bound for $(q^k)_{n\in\w}$.

\textbf{Case $\gamma$ is a limit ordinal:} Let 
$r_\gamma=\bigcup_{\beta\in M\cap\gamma}r_\beta$ and let $k\in\w$ be such that $\gamma=\alpha_k$. We must prove that $r_\gamma\le q^k(n)$ for all $n\in\w$, but it suffices to show that for all $\alpha\in\gamma\cap M$, we have $r_\gamma\restriction\alpha\le q^k(n)\restriction\alpha$. So, let $\alpha\in \gamma\cap M$ and let $m\in\w$ be such that $\alpha=\alpha_m$.

\textbf{Subcase $k>m$:} As before, by \ref{c de las q}, we know that for all $i\in\w$, 
$$q^m((k-m)+i)\le q^k(i)\restriction\alpha_m,$$
and by the inductive hypothesis $r_\gamma\restriction\alpha_m=r_{\alpha_m}\le q^m((k-m)+i)$, so we are done.

\textbf{Subcase $k<m$:} Again, by \ref{c de las q}, for all $i\in\w$, 
$$q^m(i)\le q^k((m-k)+i)\restriction\alpha_m,$$
and by the inductive hypothesis $r_{\gamma}\restriction\alpha_m=r_{\alpha_m}$ is a lower bound of $\langle q^m(i)\mid i\in\w\rangle$.

This completes the construction of the sequence $\langle r_\beta \mid \beta\in\alpha\cap M\rangle$. Finally, let 
$$r=\bigcup_{\beta\in\alpha\cap M}r_\beta.$$

\begin{subclaim}
    $r$ is a lower bound for $(p_n)_{n\in\w}$.
\end{subclaim}
\textit{Proof of the subclaim:} To see this, it suffices to prove that for every $\beta\in M\cap\alpha$ and every $n\in\w$, we have 
$$r\restriction\beta=r_\beta\le p_n\restriction\beta.$$
Suppose that $\beta=\alpha_m$. Then, by condition (\ref{ec importante}), we know that for all $n\ge m$ we have:
\[
r_\beta\le q^m(n-m)\le p_{n}\restriction\alpha_m.
\]
If $n<m$, then $p_m\le p_n$, and in particular, 
$$p_m\restriction\alpha_m\le p_n\restriction\alpha_m.$$
Thus, we conclude that $r_\beta\le p_n\restriction\alpha_m$ for every $n\in\w$. 
{\hfill\ensuremath{_{Subclaim}\square}}\finishclaim

\textbf{Construction of $\rho$ and the sequences:} The construction of the strategy $\rho$ and the sequences 
$$\langle s^m(n),q^m(n),\dot{b}^m(n) \mid m,n\in\w\rangle$$ 
will be such that in step $m=0$ we construct $s^0(0)$, $q^0(0)$, and $\dot{b}^0(0)$ (and $\rho(p_0)$ of course); in step $m=1$ we construct $s^0(1)$, $s^1(0)$, $q^0(1)$, $q^1(0)$, $\dot{b}^0(1)$, and $\dot{b}^1(0)$ (and $\rho(p_0,p_1)$); and in general, in step $m$ we construct:
\begin{itemize}
    \item $s^0(m), s^1(m-1), \dots, s^m(0)$,
    \item $q^0(m), q^1(m-1), \dots, q^m(0)$,
    \item $\dot{b}^0(m-1), \dot{b}^1(m), \dots, \dot{b}^m(0)$ and
    \item $\rho(p_0,\dots,p_m)$.
\end{itemize}

For $m=0$, let $p_0$ be the first play of $\I$, and set $s^0(0)=p_0$. Choose $q^0(0)\in\P_{\alpha_0}$ and $\dot{b}^0(0)$ in $M$ such that:
\begin{itemize}
    \item $q^0(0)\leq \sigma_0(p_0\restriction\alpha_0)$, and
    \item $q^0(0)\forces\text{``}\dot{b}^0(0)=\dot{\tau}_0(p_0(\alpha_0))\text{''}$.
\end{itemize}
Let $\rho(p_0)=\Bigl(q^0(0)\conc\dot{b}^0(0)\Bigr)\wedge p_0$. It is easy to see that $s^0(0)$, $q^0(0)$, $\dot{b}^0(0)$, and $\rho(p_0)$ satisfy all the conditions for $m=0$.

Now, suppose that we have completed the construction up to the $m$-th move of the game, and let $p_{m+1}$ be the $(m+1)$-st play of $\I$ in $\DG(\P_{\alpha_m}, M)$. Then the table for the game $\DG(\P_\alpha, M)$ now looks like:
$$
\begin{array}{c|c|c|c|c|c|c|c|c}
    \mathrm{I}  &  p_0  &                                     & p_1  & & \ldots & p_m & & p_{m+1}   \\
 \hline
  \mathrm{II}  &                    & \rho(p_0)  & & \rho(p_0,p_1) & \ldots &  & \rho(p_0,\dots,p_m)
\end{array}
$$
We now want to construct $\rho(p_0,\dots,p_{m+1})$ and the following sequences:
\begin{itemize}
    \item $s^0(m+1)$, $s^1(m)$, $\dots$, $s^{m+1}(0)$,
    \item $q^0(m+1)$, $q^1(m)$, $\dots$, $q^{m+1}(0)$, and
    \item $\dot{b}^0(m+1)$, $\dot{b}^1(m)$, $\dots$, $\dot{b}^{m+1}(0)$.
\end{itemize}

For this, enumerate $\{\alpha_0,\dots,\alpha_m,\alpha_{m+1}\}$ as $\{\gamma_0,\dots,\gamma_{m+1}\}$ in decreasing order, that is:
$$\gamma_{m+1}<\gamma_m<\dots<\gamma_0.$$
For every $n\leq m+1$, let $i_n\leq m+1$ be the unique index such that $\gamma_n=\alpha_{i_n}$.

Then, set 
$$s^{i_0}((m+1)-i_0):=p_{m+1},$$ 
and, in the game $\DG(\P_{\gamma_0}, M)$, choose $q^{i_0}((m+1)-i_0)$ and $\dot{b}^{i_0}((m+1)-i_0)$ in $M$ such that:
\begin{itemize}
    \item $q^{i_0}((m+1)-i_0)\leq \sigma_{i_0}\bigl(s^{i_0}(0)\restriction\gamma_0,\dots,s^{i_0}((m+1)-i_0)\restriction\gamma_0\bigr)$, and
    \item $q^{i_0}((m+1)-i_0)\forces\text{``}\dot{b}^{i_0}((m+1)-i_0)=\dot{\tau}_{i_0}\bigl(s^{i_0}(0)(\gamma_0),\dots,s^{i_0}((m+1)-i_0)(\gamma_0)\bigr)\text{''}$.
\end{itemize}
Thus, the table for the game $\DG(\P_{\alpha_{i_0}}, M)$ now appears as:
$$
\begin{array}{c|c|c|c|c|c|c|c}
    \mathrm{I}  &  s^{i_0}(0)\restriction\alpha_{i_0}  &                                     & \ldots & s^{i_0}(m-i_0)\restriction\alpha_{i_0} & & s^{i_0}((m+1)-i_0)\restriction\alpha_{i_0}  \\
 \hline
  \mathrm{II}  &                    & q^{i_0}(0)  & \ldots &  & q^{i_0}(m-i_0) & & q^{i_0}((m+1)-i_0)
\end{array}
$$

In general, if $q^{i_n}({(m+1)-i_n})$, $s^{i_n}({(m+1)-i_n})$ and $\dot{b}^{i_n}({(m+1)-i_n})$ have been constructed for some $n\leq m$, then we do what follows:

We know that $\gamma_{n+1}=\alpha_{i_{n+1}}$ for some $i_{n+1}\leq m+1$, then we let 
$$s^{i_{n+1}}({(m+1)-i_{n+1}}):=\left(q^{i_{n}}({(m+1)-i_n})\conc \dot{b}^{i_n}((m+1)-i_n)\right)\wedge\left(s^{i_n}({(m+1)-i_n}))\right)$$ 
and we let $q^{i_{n+1}}({(m+1)-i_{n+1}})$ and $\dot{b}^{i_{n+1}}({(m+1)-i_{n+1}})$ be such that:
\begin{itemize}
    \item $q^{i_{n+1}}({(m+1)-i_{n+1}})\leq \sigma_{i_{n+1}}(s^{i_{n+1}}(0)\restriction\gamma_{n+1},\dots,s^{i_{n+1}}({(m+1)-i_{n+1}})\restriction\gamma_{n+1})$ and
    \item $q^{i_{n+1}}({(m+1)-i_{n+1}})\forces``\dot{b}^{i_{n+1}}({(m+1)-i_{n+1}})=\dot{\tau}_{i_{n+1}}(s^{i_{n+1}}(0)(\gamma_{n+1}),\dots,s^{i_{n+1}}({(m+1)-i_{n+1}})(\gamma_{n+1}))"$.
\end{itemize}
So the table of the game $\DG(\P_{\alpha_{i_{n+1}},M})$ now looks like:
\[
\begin{array}{c|c|c|c|c|c}
    \mathrm{I}  &  s^{i_{n+1}}(0)\restriction\alpha_{i_{n+1}} & & \ldots    & s^{i_{n+1}}({(m+1)-i_{n+1}})\restriction\alpha_{i_{n+1}}  \\
 \hline
  \mathrm{II}  &                    & q^{i_{n+1}}(0)  &\ldots & &  q^{i_{n+1}}({(m+1)-i_{n+1}})
\end{array}
\]
Finally, let be $\rho(p_0,\dots,p_{m+1})$ as follows:
$$\left(q^{i_{m+1}}({(m+1)-i_{m+1}})\conc \dot{b}^{i_{m+1}}({(m+1)-i_{m+1}})\right) \wedge\dots\wedge\left(q^{i_0}({(m+1)-i_0})\conc\dot{b}^{i_0}({(m+1)-i_0})\right)\wedge p_{m+1}$$
which is equal to:
$$\left(q^{i_{m+1}}({(m+1)-i_{m+1}})\conc \dot{b}^{i_{m+1}}({(m+1)-i_{m+1}})\right)\wedge\left(s^{i_{m+1}}((m+1)-i_{m+1}\right).$$
It is not difficult to see that all the conditions are now satisfied for $m+1$. This finishes the construction of the strategy $\rho$ and the proof of the theorem.
\end{proof}

Finally, by Proposition \ref{P is statclosedinmodels}, Theorem \ref{iteracion de scm es scm} and Theorem \ref{statclosedinmodels forcings no mata estacionarios de arboles} we have that the forcing notion of Theorem \ref{diamante del complemento de S y no club en S} is \sic, so if $T$ is an $\w_1$-tree (in $V$) such that $T\restriction S\not\in NS^T$ and $T\restriction(\w_1\setminus S)\in NS^T$, $\P$ preserves these properties, that is:
$$\P\forces\text{``}(T\restriction S\not\in NS^T)\land(T\restriction (\w_1\setminus S)\in NS^T)\text{''}.$$
Thus, in $V[G]$ we have $(\diam(\w_1\setminus S))\land(\neg\diam_T)$, in particular:
$$V[G]\models \diam\wedge\left(\exists T\text{ nonspecial }\w_1\text{-tree}(  \neg\diamT\right)),$$
which proves Theorem \ref{conjectura diamante mas debil que diamante T}.

\section{Open questions}

Of course, one of the most important questions is what we can say about the diamond principle in the classical examples of nonspecial $\omega_1$-trees without cofinal branches, namely $\sigma\mathbb{Q}$ and the tree $T(S)$:

\begin{questionn}
    Does $\sigma(\mathbb{Q})$ admit a nice successor partition? More generally, does $\diam_{\sigma\mathbb{Q}}$ hold?
\end{questionn}

\begin{questionn}
    Given a bistationary $S\subseteq\omega_1$, does $\diam_{T(S)}$ hold?
\end{questionn}

On a different note, Theorem \ref{Lipschitz map and katetov} shows that a necessary condition for obtaining $\diam_T+\neg\diam_S$ is that $S\not\leq T$. A natural question is whether sufficient conditions can be given:

\begin{questionn}
    Given $\omega_1$-trees $S$ and $T$, are there sufficient conditions on $S$ and $T$ that determine when one can force $\diam_T$ together with $\neg\diam_S$?
\end{questionn}

By part \ref{2 de restrictionclubnotstationary} of Theorem \ref{restrictionclubnotstationary}, we have $NS_{\w_1}\subseteq NS_T$ for every tree $T$ of height $\w_1$. On the other hand, for every normal ideal $\mathcal{I}$ on $\w_1$ we have $NS_{\w_1}\subseteq\mathcal{I}$. A natural question, therefore, is the following:

\begin{questionn}
    Given a normal ideal $\mathcal{I}$ on $\w_1$, is there a tree $T$ of height $\w_1$ such that, for every $X\subseteq\w_1$,
    \[
    X\in\mathcal{I}\quad\text{if and only if}\quad T\restriction X\in NS^T\;?
    \]
\end{questionn}





\subsection*{Acknowledgments}
The authors are deeply grateful to Assaf Rinot and Stevo Todorcevic for their invaluable guidance, for generously taking the time to answer our questions, for suggesting directions for the investigation, and for reading earlier drafts of the manuscript. Their insightful comments significantly enhanced the content of this work. We also wish to thank Carlos Martínez Ranero for his careful reading of a previous draft; his detailed observations helped us correct an error in Lemma \ref{special bijection 2} and substantially improve the manuscript.

\bibliography{ref}{}
\bibliographystyle{plain}

\end{document}